\documentclass[onecolumn]{autart}
\usepackage{}    
\usepackage[colorlinks,linkcolor=black,anchorcolor=red,citecolor=red]{hyperref}%

\usepackage{graphicx}           
\usepackage{dcolumn}            
\usepackage{bm}                 
\usepackage{bbm}
\usepackage{graphics}           
\usepackage{epsfig}             
\usepackage{times}              
\usepackage[fleqn]{amsmath}
\usepackage{amssymb}
\usepackage{extarrows}
\usepackage{amsfonts}
\usepackage{mathrsfs}
\usepackage{fancyhdr}
\usepackage{color}
\usepackage{subfigure}
\usepackage{enumerate}
\usepackage{epstopdf}
\usepackage{subfiles}
  \graphicspath{{eps/}}
\usepackage [latin1]{inputenc} %
\usepackage{siunitx}

\usepackage{caption}
\usepackage{diagbox}

\allowdisplaybreaks[4] 

\newtheorem{theorem}{Theorem}[section] %
\newtheorem{lemma}{Lemma}[section]

\newtheorem{proposition}[theorem]{Proposition}

\newenvironment{pf3}{{\it Proof: \enspace}}{\hfill $\blacksquare$\par}

\newenvironment{pf-ISS-target}{{\it Proof  of Proposition~\ref{ISS-target}: \enspace}}{\hfill $\blacksquare$\par}
\newenvironment{pf-ISS-observer-target}{{\it Proof  of Proposition~\ref{ISS-observer-target}: \enspace}}{\hfill $\blacksquare$\par}

\begin{document}

\begin{frontmatter}

\title{Input-to-State Stabilization of $1$-D  Parabolic PDEs under  Output Feedback Control}

\author{Yongchun Bi$^{1}$}\ead{yongchun.bi@my.swjtu.edu.cn},
\author{Jun Zheng$^{1,2}$}\ead{zhengjun2014@aliyun.com},
\author{Guchuan Zhu$^{2}$}\ead{guchuan.zhu@polymtl.ca}

\address{$^{1}$School of Mathematics, Southwest Jiaotong University,
	Chengdu 611756, Sichuan, China\\
$^{2}${Department of Electrical Engineering, Polytechnique Montr\'{e}al, P.O. Box 6079, Station Centre-Ville, Montreal, QC, Canada H3T 1J4}
        }

\begin{abstract}
This paper addresses the problem of input-to-state stabilization for a class of  parabolic  equations  with time-varying coefficients, as well as  Dirichlet and Robin boundary disturbances.
By using time-invariant kernel functions, which can reduce the complexity in control design and implementation, an observer-based output feedback controller is designed via  backstepping. By using the generalized Lyapunov method, which can be used to handle Dirichlet boundary terms, the input-to-state stability of the closed-loop system under output feedback control, as well as the state estimation error system, is established in the spatial $L^\infty$-norm. Numerical simulations are conducted to confirm the theoretical results and to illustrate the effectiveness of the proposed control scheme.
\end{abstract}
\end{frontmatter}
 \section{Introduction}

 \label{sec:introduction}

In   the past few decades, backstepping is proven to be  a powerful tool for  designing  state feedback controllers, as well as output feedback controllers, in the  problem of stabilization for partial differential equations (PDEs); see, e.g.~\cite{Krstic2008,Krstic2008book}.
The main idea of backstepping is to apply the Volterra transformation to transfer the original system with destabilizing terms into a  target system, in which the  destabilizing terms  are
eliminated,  {and hence} a boundary stabilizing control  {can be obtained} from the Volterra transformation.
{The standard backstepping control design suggested in the existing literature is first to select a   target system,}   and then {to determine the corresponding} kernel function in the  Volterra {transformation. This procedure}  can make   the application  of backstepping  {in a systematic manner}, bringing convenience in solving the stabilization problem  for PDEs.
Especially,   the target system  is often chosen to be PDEs with constant coefficients {so} that the stability can be easily assessed   by, e.g., using the Lyapunov method.

However, as indicated in~\cite{chen2022arxiv}, {selecting} a target system and a kernel function  in sequence   brings    difficulties to the control design and implementation for PDEs that have  time-varying coefficients.
Indeed, when the original system has time-varying coefficients, if the target system is chosen to be with constant coefficients,  the kernel function  {will in general depend}  on the time variable (see, e.g., \cite{Karafyllis2021SIAM,Meurer2009,si2018boundary,Smyshlyaev2005AUT}),
which {leads to computation burdens and complexities in controller implementation} due to the fact that  the kernel  equation must be solved {on the fly of control execution.}

In the context of the stabilization  problem  for PDEs with external disturbances, the method of backstepping  has been used to design
boundary state feedback controllers  that can ensure the input-to-state stability (ISS) of the closed-loop system; see, e.g.,  {\cite{Karafyllis2016,Mironchenko2019,Wang2021IMA,Zheng2019TAC,Zhang2021,wen2018}}.
{It is noted that}  most of the PDEs considered in the literature involve only distributed in-domain disturbances or Robin/Neumann boundary disturbances, {while problems} with Dirichlet boundary disturbances {are much less addressed}.
The main influencing factor for this {situation} lies  {in difficulties encountered in} stability analysis. Indeed, when using backstepping to design boundary controllers for PDEs, the disturbances involved {in} the target system keep  the same type as {those} in the original system. As is well known,   in-domain disturbances or Robin/Neumann boundary disturbances do not pose difficulties for establishing the ISS   of PDEs by using the Lyapunov method (see, e.g., \cite{Mironchenko2020,Schwenninger2020}), while, as indicated in  \cite{Mironchenko2020,Zheng2024},   {the appearance of} Dirichlet boundary disturbances {represents an intrinsic obstacle} in the application of  the Lyapunov method for {ISS analysis}.

Another notable feature of the application of backstepping is that, {when} not all the state variables  are available or measured  {for state feedback control},  it is possible to  {implement output feedback control by using backstepping state observers} \cite{Krstic2008,Krstic2008book}.
{In the absence of   disturbances,}
observer-based  {backstepping} output feedback controllers  were  designed to  stabilize  single   PDEs \cite{Krstic2008,Krstic2011,smyshlyaev2005},
and   coupled PDEs or PDE-ODEs \cite{He2017,Krstic2019PDEODE,Liu2016SCL,Orlov2017}, {just to cite a few}.  However, under the framework of  ISS  theory,
{observer}-based output feedback control {design  for PDEs} in the presence of
disturbances {has not yet been reported}.  {As the presence of external disturbances is invertible in many practical application scenarios, it is obviously}   imperative to address the stabilization problem  for PDEs with disturbances under the framework of  ISS  theory, which provides a powerful tool for describing the influence of disturbances on the stability of  control systems \cite{Sontag1989,Mironchenko2020,Mironchenko2023book,Karafyllis2018book}.

In this paper, we consider the observer-based   output feedback
control  for a class of  $1$-D parabolic PDEs with time-varying coefficients and distributed in-domain  and Dirichlet-Robin boundary disturbances    under the framework of ISS theory. Unlike the approach  usually adopted in   {the} existing literature: (i) For the backstepping control  design, we apply time-invariant kernel functions to design  a state feedback controller and a state observer for output feedback control. Thus,   the complexity in control  design and implementation can be significantly reduced. A key point of   this design procedure  is that   we first select a time-invariant kernel function and then determine the target {system   with} time-varying coefficients and Dirichlet boundary disturbances,  {which is first applied in \cite{chen2022arxiv} for backstepping state feedback control.}
(ii)  {For  stability} analysis, to  deal with Dirichlet boundary terms {in} the  target {system,} we apply the so-called generalized Lyapunov functionals, which was proposed in \cite{Zheng2024},  to assess the ISS   of the involved closed-loop systems. Overall, the main contributions of this paper are threefold:
\begin{enumerate}
	\item[(i)] An observer-based output feedback controller is designed for a class of   parabolic  PDEs with disturbances  under the   framework of ISS theory.
	\item[(ii)] {The} state and output feedback  controllers are designed by using
	time-invariant kernels functions, which {can} reduce the complexity in control  design and implementation.
	
	\item[(iii)] For the involved  parabolic PDEs in the presence of various disturbances, especially, the Dirichlet boundary disturbances, the ISS  is established  in the spatial $L^\infty$-norm  by using the generalized Lyapunov method.
\end{enumerate}

In the rest of the paper, we  introduce first some basic notations used in this paper. In Section~\ref{sec:Problem formulation}, we present the problem formulation. In Section~\ref{sec:state feedback control design}, we show how to design a state feedback controller by using time-invariant kernels and prove  the first main result on the ISS of the closed-loop system by using the generalized Lyapunov method.  In Section~\ref{sec:observer design}, we design an observer  and   present the second main result on the ISS of the {state estimation} error  system.  The output feedback {control   and}   the third main result on the ISS of {the  closed-loop system under output feedback control},  as well as its proof,   are presented in Section~\ref{Sec: output-control-design}.
Numerical {simulations are conducted} in Section~\ref{sec:numerical results} to illustrate the effectiveness of the proposed control scheme. Some conclusions are given in Section~\ref{conclusion}.

	\textit{Notation:} Let  {${\mathbb{N}_0}:=\{0,1,2,...\}$, ${\mathbb{N}}:=\mathbb{N}_0\setminus\{0\}$, $\mathbb{R}:=(-\infty,+\infty)$, $\mathbb{R}_{\geq 0}:=[0,+\infty)$, $\mathbb{R}_{> 0}:=(0,+\infty)$, and $\mathbb{R}_{\leq 0}:=(-\infty,0]$.}

	For a domain $\Omega$  (either open or {closed}) in $\mathbb{R}^1$ or $\mathbb{R}^2$, let $C\left({\Omega}\right):=\{v: {\Omega} \rightarrow \mathbb{R}|~ v$ is continuous on $\Omega\}$.
	For $i\in\mathbb{N}$, $C^i\left({\Omega}\right):=\{v: {\Omega} \rightarrow\mathbb{R}|~ v$ has continuous derivatives up to order $i$ on ${\Omega}$\}.  Let $C^{2,1}(\Omega)$$:=$$\{v:$ $\Omega $$\rightarrow$$\mathbb{R}|~$$v$, $v_x$, $v_{xx}$, $v_t$$\in$ $C(\Omega)\}$.

	For $T\in\mathbb{R}_{>0} $, let $Q_T := (0,1) \times  (0,T)$ and $\overline{Q}_T$$:=$$[0,1]$$\times$$[0,T]$.  Let $Q_\infty$$:=(0,1)\times {\mathbb{R}_{>0}}$ and $\overline{Q}_\infty:=[0,1]\times {\mathbb{R}_{\geq 0}}$.
	For $v: Z \rightarrow \mathbb{R}$ with $Z \subset \overline{Q}_{\infty}$, the notation $v[t]$ {(or $v[{z}]$)} denotes the profile at certain $t \in \mathbb{R}_{\geq 0}$ (or ${z} \in[0,1]$), i.e., $v[t]({z})=v({z}, t)$ (or $v[{z}](t)=v({z}, t)$).
	
	The norms of a function $v$ in  the Lebesgue spaces $L^p(0,1)$ with  $p\in [1,+\infty)$  and $L^{\infty}(0,1)$  are defined by $\|v\|_{L^p(0,1)}:=\left(\int_0^1|v(x)|\text{d}x\right)^{\frac{1}{p}}$ and $\|v\|_{L^{\infty}(0,1)}:= \operatorname{ess \sup }_{x \in(0,1)} |v(x)|   $, respectively.
	For $T\in\mathbb{R}_{>0} $, the norm of a function $v$ in the  Lebesgue space  $L^{\infty}((0, T) ; L^{\infty}(0,1))$  is defined by
	$\|v\|_{L^\infty((0, T) ; L^{\infty}(0,1))}:=\operatorname{ess \sup }_{t \in(0, T)}\|v[t]\|_{L^{\infty}(0,1)}$.

	Let $\mathcal {K} := \{\gamma:\mathbb{R}_{\geq 0} \rightarrow \mathbb{R}_{\geq 0}|$$\gamma(0)${$=$}{$0$}\}, $\gamma$ is continuous, strictly increasing\}, $\mathcal {L}$$:=$$\{\gamma:$$\mathbb{R}_{\geq 0}$$\rightarrow$$\mathbb{R}_{\geq 0}|$ $\gamma$ is continuous, strictly decreasing, {$\lim_{s\rightarrow\infty}\gamma(s)$}{$=$}{$0$}\}, $\mathcal {KL} := \{\beta:$$\mathbb{R}_{\geq 0}$$\times$$\mathbb{R}_{\geq 0}$ $\rightarrow$$\mathbb{R}_{\geq 0}|$${\beta}$ is continuous, $\beta[t]\in\mathcal {K}$, $\forall t \in \mathbb{R}_{\geq 0}$; $\beta[s] \in \mathcal {L},\forall  s \in {\mathbb{R}_{> 0}}\}$.

\section{Problem Statement}\label{sec:Problem formulation}
In this paper, we address the problem of stabilization of one-dimensional linear parabolic equation:
\begin{subequations}\label{original system}
	\begin{align}
		u_t(x,t)=&u_{xx}(x,t)+\lambda(t)u(x,t)+f(x,t), (x,t)\in Q_\infty,\label{1a}\\
		u_x(0,t)=&qu(0,t)+d_0(t), t\in \mathbb{R}_{>0},\label{1b}\\
		u(1,t)=&U(t)+d_1(t), t\in \mathbb{R}_{>0},\\
		u(x,0)=&u_0(x), x\in(0,1),
	\end{align}
\end{subequations}
where   $\lambda(t) $ is the reaction coefficient, which is time-varying,  $f(x,t)$ denotes in-domain disturbance, $d_0(t)$ and $d_1(t)$  {denote} Robin and Dirichlet boundary disturbance, respectively,   $q $ is a constant,  $u_0(x)$ is the initial datum, and    $U(t)$ is a boundary input  that needs to be designed.

The main {objective} of this work {is}
to design  an observer-based output   feedback  control law  $U(t)$ via    time-invariant kernel functions such that  both  the {state estimation} error system and {the   closed-loop system} are ISS  with respect to (w.r.t.) external {disturbances.}

Throughout this paper, we always make the following assumptions:
\begin{subequations}
	\begin{align}
		q\in &\mathbb{R}_{>0},
		\lambda\in  C(\mathbb{R}_{\ge0}),\mathop{\text{sup}}\limits_{t\in\mathbb{R}_{\ge 0}}\lambda(t)<+\infty,\label{assumpations}\\
		f \in&C \left(\overline{Q}_{\infty}\right), d_0 , d_1 \in C (\mathbb{R}_{\geq 0}) ,\\
		u_0\in& {C^{2,1}(0,1)}\cap C([0,1]).
	\end{align}
\end{subequations}
Moreover, the initial datum $u_0 $ and boundary disturbances $d_0,d_1$ satisfy  certain compatibility conditions, which will be specified later in each section.


\section{State Feedback {Control  and} ISS of the Closed-loop System}\label{sec:state feedback control design}
Before designing an {observer for output} feedback {control}, we show {first} that there {exists a  state feedback control}  that can ensure the ISS in the  $L^{\infty}$-norm {of system~\eqref{original system}}.
\subsection{{State Feedback Control Design}}
\subsubsection{Selection of Kernel Functions}
To adopt time-invariant kernel functions in backstepping control design,  following~\cite{chen2022arxiv}, we first {select  a kernel} function, and then {select a target} {system  with time-varying} coefficient. More specifically, by  the assumption~\eqref{assumpations}, we  take any positive constant $c_0$ such that
\begin{align}\label{c0}
	c_0>\mathop{\text{sup}}\limits_{t\in\mathbb{R}_{\ge 0}}\lambda(t).
\end{align}
Denote
\begin{align}\label{c lambda}
	c(t):=c_0-\lambda(t) \text{ and } \underline{c}:=\inf_{t\in {\mathbb R_{\ge 0}}} c(t)>0.
\end{align}

Let $D{:=}\{(x,{z})\in \mathbb{R}^2|~0\leq {z} \leq x\leq 1\}$. We consider the solution $k$ {to} the following equation defined over $D$:
\begin{subequations}\label{original kernel equation}
	\begin{align}
		k_{xx}(x,{z})-k_{{zz}}(x,{z})=&c_0k(x,{z}),\\
		\frac{\text{d}}{\text{d}x}(k(x,x))=&-\frac{1}{2}c_0,\\
		k_{{z}}(x,0)=&qk(x,0),\\
		k(0,0)=&0,
	\end{align}
\end{subequations}
where $\frac{\text{d}}{\text{d}x}\left(k(x,x)\right) :=k_x(x,{z})|_{{z}=x}+k_{{z}}(x,{z})|_{{z}=x}$. Indeed, the existence and regularity of the solution $k$ to equation~\eqref{original kernel equation} are guaranteed by the following Lemma, whose proof can be found in \cite{smyshiyaev2004}.
\begin{lemma}\label{proposition kernel}
	Equation~\eqref{original kernel equation} has a unique solution $k\in C^{2}\left(D\right)$, which can be expressed as
	\begin{align}\label{specific-k}
		 k(x, {z})=&-c_0 x \frac{\mathcal{I}_1\left(\sqrt{c_0\left(x^2-{z}^2\right)}\right)}{\sqrt{c_0\left(x^2-{z}^2\right)}}\notag\\
&+\frac{q c_0}{\sqrt{c_0+q^2}} \int_0^{x-{z}} e^{\frac{-q \tau} {2}}\mathcal{I}_0 \left(\sqrt{c_0(x+{z})(x-{z}-\tau)}\right)\sinh \left(\frac{\sqrt{c_0+q^2}}{2} \tau\right) {\rm d} \tau,
	\end{align}
	where $\mathcal{I}_0 $ and $\mathcal{I}_1 $ are the modified Bessel functions of order zero and order one, respectively {(see, e.g.,~\cite[pp. 174-175]{Krstic2008book})}.
\end{lemma}

\subsubsection{Selection {of Target} System}

Applying the Volterra transformation
\begin{align}\label{original transfer}
	w(x,t):=u(x,t)-\int_{0}^{x}k(x,{z})u({z},t)\text{d}{z},
\end{align}
we obtain the following target system   involving a time-varying reaction coefficient and also   Dirichlet and Robin boundary disturbances:
\begin{subequations}\label{target system}
	\begin{align}\label{target1}
		w_t(x,t)=&w_{xx}(x,t)-c(t)w(x,t)+\psi(x,t), &&(x,t)\in Q_{\infty},\\
		w_x(0,t)=&qw(0,t)+d_0(t), &&t\in \mathbb{R}_{>0},\\
		w(1,t)= &d_1(t), &&t\in \mathbb{R}_{>0},\\
		w(x,0)=&w_0(x), &&x\in(0,1),
	\end{align}
\end{subequations}
where
\begin{align}
	\psi(x,t){:=} &f(x,t)-\int_{0}^{x}k(x,{z}) f({z},t)\text{d}{z}+k(x,0) d_0(t),\label{def-psi}\\
	w_0(x){:=}&u_0(x)-\int_{0}^{x}k(x,{z})u_0({z})\text{d}{z}.\notag
\end{align}
		
Note that the transformation~\eqref{original transfer} is invertible, and the equivalence between system~\eqref{original system} and   system~\eqref{target system}  is guaranteed by the following lemma, whose  proof   can be proceeded in a standard way as in, e.g., \cite{Coron2017}  and hence, is omitted.

\begin{lemma}\label{lem.1'}
	Let $l$ be the unique solution of
	\begin{subequations}\label{inverse kernel equation}
		\begin{align}
			l_{xx}(x,{z})-l_{{zz}}(x,{z})=&-c_0l(x,{z}),\label{inverse equation1}\\
			\frac{\text{d}}{\text{d}x}(l(x,x))=&-\frac{1}{2}c_0,\label{inverse equation2}\\
			l_{{z}}(x,0)=&ql(x,0),\\
			l(0,0)=&0.
		\end{align}
	\end{subequations}
	Then,
	the inverse transformation of~\eqref{original transfer}  is given by
	\begin{align}\label{inverse transfer conclusion}
		u(x,t):=w(x,t)+\int_0^x l(x, {z}) w({z},t) \mathrm{d} {z}.
	\end{align}
	Therefore,    system~\eqref{original system} is equivalent to  system~\eqref{target system}.
\end{lemma}

	\subsubsection{State Feedback Controller}
	By using the function $k$ satisfying \eqref{original kernel equation} and in view of the transformation~\eqref{original transfer}, we   define the state feedback control law as
	\begin{align}\label{control law}
		U(t):=\int_{0}^{1}k(1,{z})u({z},t)\text{d}{z}.
	\end{align}
	\subsection{ISS Assessment of the Closed-loop System}
	{In this section, we assess the ISS of the closed-loop {system} under the state feedback control law \eqref{control law}}.
	\subsubsection{ISS of the Target System}
	{To  ensure} the existence of a classical solution to the original system~\eqref{original system} and the  target system~\eqref{target system},   we assume further that the initial datum $u_0$ and the boundary disturbances $d_0,d_1$ satisfy  the  compatibility conditions:
	\begin{align*}
		u_{0x}(0)= qu_0(0)+d_0(0),
		u_0(1)= \int_{0}^{1}k(1,{z})u_0({z})\text{d}{z}+d_1(0).
	\end{align*}
	The  classical PDE theory (see, e.g., \cite[Chapter  \uppercase\expandafter{\romannumeral4}]{Ladyzenskaja1968}) ensures that  system~\eqref{target system}  admits a unique solution $w\in C^{2,1}(Q_T)\cap C\left(\overline{Q}_T\right)$ and hence, system~\eqref{original system}  admits a unique solution $u\in C^{2,1}(Q_T)\cap C\left(\overline{Q}_T\right)$ for any $T\in \mathbb{R}_{>0}$.
	In addition, for the target system~\eqref{target system}, we can prove the following result, which indicates the ISS w.r.t. external disturbances.
	\begin{proposition}\label{ISS-target} The target system~\eqref{target system} is ISS  in the  $L^\infty$-norm w.r.t.  $f,d_0,d_1$, having the following estimate for all $T \in \mathbb{R}_{> 0}$:
		{\begin{align*}
				\left\|{w}[T]\right\|_{L^{\infty}(0,1)}
				\leq c_1\left(e^{-\sigma T}\left\|{w}_0\right\|_{L^\infty(0,1)}+ \left\|{d}_0\right\|_{L^{\infty}(0, T)}\right.\left.+\left\|{d}_1\right\|_{L^{\infty}(0, T)}+\left\|{f}\right\|_{L^{\infty}\left((0,T);L^\infty(0,1)\right)}\right),
			\end{align*}
			where $\sigma$ is an arbitrary constant satisfying $\sigma\in (0,\underline{c})$,  and $c_1$  is a positive constant depending only on $q,\underline c,\sigma$, and $\max_{(x,{z})\in D}|k(x,{z})|$.}
	\end{proposition}
	
	{The proof of Proposition~\ref{ISS-target} is given in Section~\ref{Sec: ISS-assessment-target}.}

\subsubsection{ISS of the Closed-loop System}
By virtue of Lemma~\ref{lem.1'} and Proposition~\ref{ISS-target}, we state  the first main result, which indicates the ISS of {the  closed-loop system~\eqref{original system} under} state feedback control.
\begin{theorem}\label{original main result}
	Under the state feedback control law~\eqref{control law}, the closed-loop system~\eqref{original system}  is ISS in the  $L^\infty$-norm w.r.t.  $f,d_0,d_1$, having the following estimate  for all $T \in \mathbb{R}_{> 0}$:
	\begin{align*}
		\left\|{u}[T]\right\|_{L^{\infty}(0,1)}
		\leq {c_{2}}\left(e^{-\sigma T}\left\|{u}_0\right\|_{L^\infty(0,1)}+\left\|{d}_0\right\|_{L^{\infty}(0, T)}+\left\|{d}_1\right\|_{L^{\infty}(0, T)}+\left\|{f}\right\|_{L^{\infty}\left((0,T);L^\infty(0,1)\right)}\right),
	\end{align*}
	where $\sigma$ is an arbitrary constant satisfying $\sigma\in (0,\underline{c})$, and {$c_{2}$}  is a positive constant  depending only on $q$,   $\underline{c}$, $\sigma$, $\max _{(x,{z})\in D}|k(x, {z})|$, {and $\max _{(x,{z})\in D}|l(x, {z})|$}.
\end{theorem}

{The proof of Theorem~\ref{original main result} is standard and hence, is omitted.}


\subsubsection{Proof of Proposition~\ref{ISS-target}} \label{Sec: ISS-assessment-target}
{In this section, we apply  the generalized Lyapunov method, {which is} proposed in \cite{Zheng2024},  to prove the ISS of the target system~\eqref{target system}.}

\begin{pf-ISS-target} We proceed with the proof in four steps.
	
	\textbf{Step 1:} define truncation functions.
	For an  arbitrary constant $r>1$, define
	\begin{align}\label{gG}
		g(\theta ):= \begin{cases}\theta^r, & \theta  \geq 0 \\
			0, & \theta <0\end{cases}\ \text{and}\
		G(\theta ):=\int_0^\theta  g(\tau) \mathrm{d} \tau,
	\end{align}
	which satisfy
	\begin{subequations}\label{properties}
		\begin{align}
			g(\theta )  \geq& 0, g^{\prime}(\theta ) \geq 0, G(\theta ) \geq 0,  \forall \theta  \in \mathbb{R},\label{gG1} \\
			g(\theta ) =&G(\theta )=0, \forall \theta \in \mathbb{R}_{\leq 0}.\label{gG2}
		\end{align}
	\end{subequations}
	
	\textbf{Step 2:} prove the upper bound of ${w}$.
	For any ${\sigma}\in (0,\underline {c})$, let
	\begin{align*}
		{{v}(x, t):=}e^{{\sigma} t} {w}(x, t),{v}_0(x):={w}_0(x),
		{\check{d}_0(t):=} e^{{\sigma} t} {d}_0(t), \check{d}_1(t):= e^{{\sigma} t} {d}_1(t),
		\check{\psi}(x, t):= e^{{\sigma} t} {\psi}(x, t).
	\end{align*}
	
	For any $ T\in\mathbb{R}_{>0}$, 
%
	let $\check{\mathcal{D}}:= \max \left\{\left\|{v}_0\right\|_{L^\infty(0,1)},\frac{1}{{q}}\left\|\check {d}_0\right\|_{L^{\infty}(0, T)},\left\|\check{d}_1\right\|_{L^{\infty}(0, T)}, \frac {1}{\underline{c}-{\sigma}}\left\|\check{\psi}\right\|_{L^{\infty}\left((0,T);L^\infty(0,1)\right)}\right\}.$
	Note that
	\begin{align}\label{v le hat k}
		{v}(1,t)\leq\check{\mathcal{D}},\forall t\in(0,T).
	\end{align}

	By     integrating by parts and using the properties of $g$ and $G$,  we get
	\begin{align}\label{widetildev}
		\frac{\text{d}}{\text{d} t} \int_0^1 G\left({v}-\check{\mathcal{D}}\right) \text{d} x
		= &g\left({v}(1, t)-\check{\mathcal{D}}\right) {v}_x(1, t) -g\left({v}(0, t)-\check{\mathcal{D}}\right)\left({q}{v}(0, t)+\check{d}_0(t)\right) -\int_0^1 g^{\prime}\left({v}-\check{\mathcal{D}}\right) {v}_x^2  \text{d} x\notag\\
		&+ \int_0^1 g\left({v}-\check{\mathcal{D}}\right)\left(({\sigma}- c(t)) {v} +\check{\psi} \right) \text{d} x \notag\\
		\leq&0,\forall t\in (0,T).
	\end{align}
	Therefore,   it holds that
	\begin{align*}
		\int_0^1G\left({v}(x,T)-\check{\mathcal{D}}\right)\text{d}x \leq \int_0^1G\left({v}(x,0)-\check{\mathcal{D}}\right)\text{d}x,
	\end{align*}
	which, along with the continuity of ${v}$, implies that
	\begin{align*}
		{v}(x,T)\leq\check{\mathcal{D}},\forall x\in(0,1).
	\end{align*}
	Then, we have
	\begin{align}\label{upper}
		 {w}(x, T)=e^{-{\sigma} T} {v}(x, T)
		\leq e^{-\sigma T}\left\|{w}_0\right\|_{L^{\infty}(0,1)}+\frac{1}{{q}}\left\|{d}_0\right\|_{L^{\infty}(0, T)}  + \left\|{d}_1\right\|_{L^{\infty}(0, T)}+\frac{1}{\underline{c}-{\sigma}}\left\|{\psi}\right\|_{L^{\infty}\left((0,T);L^\infty(0,1)\right)}.
	\end{align}

	\textbf{Step 3:} prove the lower bound of {${w}$}.
	For the same ${\sigma}\in (0,\underline {c})$  as in  Step 2, let
	\begin{align*}
		{\overline{v}(x, t):=}-e^{{\sigma} t}{w}(x, t),\overline{v}_0(x):=-{w}_0(x),{\overline{d}_0(t):=} -e^{{\sigma} t} {d}_0(t), \overline{d}_1(t):= -e^{{\sigma} t} {d}_1(t),\overline{\psi}(x, t):= -e^{{\sigma} t} {\psi}(x, t).
	\end{align*}
 	Let $g$, $G$ be defined by~\eqref{gG}, and
$\overline{\mathcal{D}}:=\max \left\{\left\|\overline{v}_0\right\|_{L^{\infty}(0,1)}, \frac{1}{{q}}\left\|\overline{d}_0\right\|_{L^{\infty}(0, T)},\left\|\overline{d}_1\right\|_{L^{\infty}(0, T)},\frac{1}{\underline{c}-{\sigma}}\left\|\overline{\psi}\right\|_{L^{\infty}((0,T);L^\infty(0,1))}\right\}.$
	
	Analogous to Step 2, for any $t\in(0,T)$, by deriving
$\frac{\text{d}}{\text{d} t} \int_0^1 G\left(\overline{v}-\overline{\mathcal{D}}\right) \text{d} x \leq 0,$
	we get:
	\begin{align*}
		\int_0^1G\left(\overline{v}(x, T)-\overline{\mathcal{D}}\right) \text{d} x \leq \int_0^1 G\left(\overline{v}(x, 0)-\overline{\mathcal{D}}\right) \text{d} x.
	\end{align*}
	Therefore, it holds that
	\begin{align*}
		\overline{v}(x, T) \leq \overline{\mathcal{D}}, \forall x \in(0,1) ,
	\end{align*}
	%
	%
	%
	%
	which implies that
	\begin{align}\label{lower}
		{-{w}(x, T)}
		\leq e^{-{\sigma} T}\left\|{w}_0\right\|_{L^\infty(0,1)}+\frac{1}{q}\left\|{d}_0\right\|_{L^{\infty}(0, T)}  +\left\|{d}_1\right\|_{L^{\infty}(0 , T)}+\frac{1}{\underline{c}-{\sigma}}\left\|{\psi}\right\|_{L^{\infty}((0, T);L^\infty(0,1))}.
	\end{align}
		
	\textbf{Step 4:} conclusion. {By the definition of $\psi$ (see~\eqref{def-psi}),} we have
	\begin{align}\label{defenition psi}
		{\left\|{\psi}\right\|_{L^{\infty}((0, T);L^\infty(0,1))}}\le \overline{k}\left\|{d}_0\right\|_{L^{\infty}(0, T)}+\left(\overline{k}+1\right)\left\|{f}\right\|_{L^{\infty}((0, T);L^\infty(0,1))},
	\end{align}
	{where $\overline{k}:=\max _{(x,{z})\in D}|k(x, {z})|$.}
	
	{By \eqref{upper}, \eqref{lower}, and \eqref{defenition psi},  we deduce that the target system~\eqref{target system} is ISS in the  $L^{\infty}$-norm w.r.t.  $f, d_0$, and $ d_1$.}
\end{pf-ISS-target}


\section{Observer  Design and ISS of the State Estimation Error System}\label{sec:observer design}
In this section, we design an observer, {which ensures} the ISS in the $L^\infty$-norm of the {state estimation} error system.
\subsection{Observer  Design}
We  design an observer in the anti-collocated {setup where} the sensor and the actuator are placed at the opposite ends.
More specifically, for system~\eqref{original system} with Dirichlet actuation, we assume that the only available measurement of the system is {located} at $x=0$, the opposite end to {the} actuation {located} at $x=1${. Thus,}  $u(0,t)$ is the {system} output. Note that in practice  the output may also be affected by external disturbances. Thus,  {denoting the measurement disturbance} by ${d_m(t)}$ and considering the output
\begin{align}
	{y}(t):=u(0,t)+{d_m(t)},\label{output}
\end{align}
we propose an observer under the following form:
\begin{subequations}\label{observer system}
	\begin{align}
		\hat{u}_t(x,t)=&\hat{u}_{xx}(x,t)+\lambda(t)\hat{u}(x,t) + p(x)\left({y}(t)-\hat{u}(0,t)\right), (x,t)\in Q_\infty,\label{observer system1a}\\
		\hat{u}_x(0,t)=&q{\hat u}(0,t) -p_{0}\left({y}(t)-\hat{u}(0,t)\right), t\in \mathbb{R}_{>0},\label{observer system2a}\\
		\hat{u}(1,t)=&U(t) , t\in \mathbb{R}_{>0},\label{observer system3a}\\
		\hat{u}(x,0)=&\hat{u}_0(x), x\in(0,1),\label{observer system4a}
	\end{align}
\end{subequations}
where  $\hat u_0(x)$ is the initial datum,  and the function $p(x)$ and the constant $p_{0}$ are observer gains to be determined.

To determine the observer gains $p(x)$ and $p_0$, we first
consider    the following equation defined over $D$:
\begin{subequations}\label{error kernel equation}
	\begin{align}
		m_{xx}(x,{z})-m_{{zz}}(x,{z})=&c_0m(x,{z}),\label{error kernel equation1}\\
		\frac{\text{d}}{\text{d}x}\left(m(x,x)\right)=&-\frac{1}{2}c_0,\label{error kernel equation2}\\
		m(1,{z})=&0,
	\end{align}
\end{subequations}
where $\frac{\text{d}}{\text{d}x}\left(m(x,x)\right) :=m_x(x,{z})|_{{z}=x}+m_{{z}}(x,{z})|_{{z}=x}$. The existence and regularity of a solution $m$ to \eqref{error kernel equation} are guaranteed by the following Lemma.
\begin{lemma}\label{kernel proposition}
	Equation~\eqref{error kernel equation} admits a unique solution $m\in C^{2}\left(D\right)$, which can be expressed as
	\begin{align*}
		m(x,{z})=c_0(1-x)\frac{\mathcal{J}_1\left(\sqrt{c_0({z}^2-x^2+2(x-{z}))}\right)}{\sqrt{c_0\left({z}^2-x^2+2(x-{z})\right)}},
	\end{align*}
	where $\mathcal{J}_1 $ is the Bessel function of order one  (see, e.g.,~\cite[pp.~173-174]{Krstic2008book}).
\end{lemma}

\begin{pf3} Let us make a change of variables
	\begin{align*}
		\overline{x}:=1-{z}, \overline{{z}}:=1-x, \overline{m}(\overline{x}, \overline{{z}}):=m(1-\overline{{z}},1-\overline{x}).
	\end{align*}
	Then, \eqref{error kernel equation} becomes
	\begin{align*} \overline{m}_{\overline{x}\overline{x}}(\overline{x},\overline{{z}})-\overline{m}_{\overline{{z}}\overline{{z}}}(\overline{x},\overline{{z}})=&-c_0\overline{m}(\overline{x},\overline{{z}}),\\
		\frac{\text{d}}{\text{d}\overline{x}}\left(\overline{m}(\overline{x},\overline{x})\right)=&\frac{1}{2}c_0,\\
		\overline{m}(\overline{x},0)=&0,
	\end{align*}
	whose solution  is expressed as (see {\cite[(4.33), p.~35]{Krstic2008book}})
	\begin{align*} \overline{m}(\overline{x},\overline{{z}})=c_0\overline{{z}}\frac{\mathcal{I}_1\left(\sqrt{-c_0(\overline{x}^2-\overline{{z}}^2)}\right)}{\sqrt{-c_0(\overline{x}^2-\overline{{z}}^2)}}=c_0\overline{{z}}\frac{\mathcal{I}_1\left(\text{i}\sqrt{c_0(\overline{x}^2-\overline{{z}}^2)}\right)}{\text{i}\sqrt{c_0(\overline{x}^2-\overline{{z}}^2)}},
	\end{align*}
	where { $\mathcal{I}_1 $ is the modified Bessel function of order one.
		Then, using the fact that $\mathcal I_1(\text{i}x)=\text{i}\mathcal J_1(x)$ (see \cite[(A.11), p.~174]{Krstic2008book}), we obtain}
	\begin{align*} \overline{m}(\overline{x},\overline{{z}})=c_0\overline{{z}}\frac{\mathcal{J}_1\left(\sqrt{c_0(\overline{x}^2-\overline{{z}}^2)}\right)}{\sqrt{c_0(\overline{x}^2-\overline{{z}}^2)}}.
	\end{align*}
\end{pf3}

 Now, based on the function $m$,   we choose the observer gains $p_0$  and $p(x)$ that  satisfy
\begin{align}
	{p_{0}} >&m(0,0)-q=\frac{c_0}{2}-q,\label{p10}\\
	p(x) =&-m_{{z}}(x,0)+{(q+p_{0})}m(x,0)+\int_{0}^{x}m(x,{z})p({z})\text{d}{z}.\label{p1}
\end{align}
The existence and regularity of a solution $p$ to the integral equation~\eqref{p1} are guaranteed by the following Lemma.

\begin{lemma}[{\cite[Theorem~1.3.1, p. 33]{corduneanu1991}}]\label{kernelpx proposition}
	Equation~\eqref{p1} admits a unique solution $p\in C(D)$.  Furthermore, $p\in C^1(D)$ due to  the fact that $m\in C^2(D)$.
\end{lemma}
\subsection{ISS of the {State Estimation} Error System}
{Let the estimation error of the observer given in \eqref{observer system} be defined by}
\begin{align}\label{observer error}
	\widetilde{u}(x,t):=u(x,t)-\hat{u}(x,t).
\end{align}
{Thus, the state estimation error system will be given by}
\begin{subequations}\label{error system}
	\begin{align}
		\widetilde{u}_t(x,t)=&\widetilde{u}_{xx}(x,t)+\lambda(t)\widetilde{u}(x,t)+\widetilde{f}(x,t)-p(x)\widetilde{u}(0,t), (x,t)\in Q_\infty,\label{error system1a}\\
		\widetilde{u}_x(0,t)=&\widetilde{d}_0(t)+{(q+p_{0})}\widetilde{u}(0,t),t\in \mathbb{R}_{>0},\\
		\widetilde{u}(1,t)=& \widetilde{d}_1(t), t\in \mathbb{R}_{>0},\label{error system1b}\\
		\widetilde{u}(x,0)=&\widetilde{u}_0(x), x\in(0,1),\label{error system1c}
	\end{align}
\end{subequations}
where $\widetilde{f}(x,t):=f(x,t)-p(x){d_m(t)}, \widetilde{d}_0(t):=d_0(t)+p_0{d_m(t)}, \widetilde{d}_1(t):=d_1(t) $, and $\widetilde u_0(x):=u_0(x)-\hat u_0(x)$.

%

By using the transformation
\begin{align}\label{transfer}
	\widetilde{w}(x,t):=\widetilde{u}(x,t)-\int_{0}^{x}m(x,{z})\widetilde{u}({z},t)\text{d}{z},
\end{align}
we obtain the following target system involving also a time-varying reaction coefficient {and    external disturbances}:
\begin{subequations}\label{error target system}
	\begin{align}
		\widetilde{w}_t(x,t)=&\widetilde{w}_{xx}(x,t)-{c}(t)\widetilde{w}(x,t)+\widetilde{\psi}(x,t), (x,t)\in Q_\infty,\label{error target system1a}\\
		\widetilde{w}_x(0,t)=&b\widetilde{w}(0,t)+\widetilde{d}_0(t),t\in \mathbb{R}_{>0},\\
		\widetilde{w}(1,t)= &\widetilde{d}_1(t), t\in \mathbb{R}_{>0},\label{error target system1b}\\
		\widetilde{w}(x,0)=&\widetilde{w}_0(x), x\in(0,1),\label{error target system1c}
	\end{align}
\end{subequations}
where
\begin{align}
	b:=&{q+p_{0}}-m(0,0)>0\ {(\text{by}\ \eqref{p10}),}\label{def-b}\\ \widetilde{\psi}(x,t):=&\widetilde{f}(x,t)-\int_{0}^{x}m(x,{z})\widetilde{f}({z},t)\text{d}{z}+m(x,0)\widetilde{d}_0(t),\label{def-tilde-psi}\\ \widetilde{w}_0(x):=&\widetilde{u}_0(x)-\int_{0}^{x}m(x,{z})\widetilde{u}_0({z})\text{d}{z}.\notag
\end{align}

 The equivalence between the state estimation error system~\eqref{error system}   and the corresponding target system~\eqref{error target system} is guaranteed by the following lemma, whose proof is standard and hence, is omitted, and the existence and regularity of a solution $n$ to the following kernel function are guaranteed by~\cite{smyshiyaev2004}.
\begin{lemma}\label{lem.2'}
	Let $n$ be the unique solution of
	\begin{subequations}\label{inverse kernel equation2}
		\begin{align}
			n_{xx}(x,{z})-n_{{zz}}(x,{z})=&-c_0n(x,{z}),\label{inverse equation3}\\
			\frac{\text{d}}{\text{d}x}(n(x,x))=&-\frac{1}{2}c_0,\label{inverse equation4}\\
			n(1,{z})=&0.
		\end{align}
	\end{subequations}
	Then, the inverse transformation of \eqref{transfer} is given by
	\begin{align}\label{inverse transfer conclusion2}
		\widetilde u(x,t):=\widetilde w(x,t)+\int_0^x n(x,{z}) \widetilde w({z},t) \mathrm{d}{z}.
	\end{align}
	Therefore,   the state estimation error system~\eqref{error system}  is equivalent to the corresponding target system~\eqref{error target system}.
\end{lemma}

To ensure the existence  of  a solution $\widetilde{u}$ to the error  system~\eqref{error system},   as well as a solution $(u,\hat u)$ to the system  coupled by \eqref{original system} and \eqref{observer system}, we assume   that the disturbance
${d_m}  \in C (\mathbb{R}_{\geq 0})$ and the initial datum $
\hat{u}_0\in  C^{2,1}(0,1)\cap C([0,1])
$
satisfy the following compatibility conditions:
\begin{align*}
	\hat{u}_{0x}(0)=& q\hat{u}_0(0)
	-p_{0}\left(u_0(0)+{d_m}(0)-\hat{u}_0(0)\right),
	\hat{u}_0(1)=  U(0).
\end{align*}

Analogous to the target  system~\eqref{target system} of the   $u$-system,  for any $T\in \mathbb{R}_{>0}$,  the target system~\eqref{error target system}   {admits a unique solution} $\widetilde w\in C^{2,1}(Q_T)\cap C\left(\overline{Q}_T\right)$ and hence, the system  coupled by the  original system~\eqref{original system} and the observer system~\eqref{observer system}  admits a unique solution $(u,\hat u)\in \left(C^{2,1}(Q_T)\cap C\left(\overline{Q}_T\right)\right)\times\left(C^{2,1}(Q_T)\cap C\left(\overline{Q}_T\right)\right) $.

{Note that {the target} system~\eqref{error target system} has the same form as  the target system~\eqref{target system} of the plant~\eqref{original system}. By virtue of Lemma~\ref{lem.2'} and Proposition~\ref{ISS-target},    we have the following theorem, which is the second main result. The proof is omitted as it is standard.}
\begin{theorem}\label{error main result} The {state estimation error} system~\eqref{error system}  is ISS in the  $L^\infty$-norm w.r.t.  $f,{d_m},d_0,d_1$,
having the following estimate  for all $T \in \mathbb{R}_{> 0}$:
\begin{align*}
	\left\|\widetilde{u}[T]\right\|_{L^{\infty}(0,1)}
	\leq    c_3 \left(e^{-\sigma T} \left\|\widetilde{u}_0\right\|_{L^\infty(0,1)}+\left\|{d_m} \right\|_{L^{\infty}(0, T)}+\left\| {d}_0\right\|_{L^{\infty}(0, T)}   +\left\| {d}_1\right\|_{L^{\infty}(0, T)}+\left\|{f}\right\|_{L^{\infty}\left((0,T);L^\infty(0,1)\right)}\right) ,
\end{align*}
where $\sigma$ is an arbitrary constant satisfying $\sigma\in (0,\underline{c})$,  and $c_3$  is a positive constant depending only on $\underline{c},  \sigma, p_{0}, {q}$,   {$\max_{x\in[0,1]}|p(x)|$,}  $\max _{(x,{z})\in D}|m(x, {z})|$, and {$\max _{(x,{z})\in D}|n(x,{z})|$}.
\end{theorem}

\section{Output Feedback {Control   and} ISS of the Closed-loop System}\label{Sec: output-control-design}
\subsection{Observer-based Output Feedback {Control}}
In this section, we {present} an observer-based output feedback {control scheme}, which ensures the ISS in the $L^\infty$-norm of the closed-loop system.

Indeed, by using the function $k$ in {closed form}~\eqref{specific-k}  and replacing $u$ with $\hat{u}$ in \eqref{control law}, we {obtain} the observer-based output feedback control {law:}
\begin{align}\label{observer control law}
U(t):=\int_{0}^{1}k(1,{z})\hat u({z},t)\text{d}{z}.
\end{align}

Note that by using the  transformation
\begin{align}\label{observer transfer}
\hat w(x,t):=\hat u(x,t)-\int_{0}^{x}k(x,{z})\hat u({z},t)\text{d}{z},
\end{align}
we can transform the {observer~\eqref{observer system}} into
\begin{subequations}\label{observer target system}
\begin{align}
	\hat w_t(x,t)=&\hat w_{xx}(x,t)-c(t)\hat w(x,t){+\breve\psi(x,t)}{+K_p(x)\widetilde{u}(0,t)}, (x,t)\in Q_\infty,\label{observer target system1a}\\
	\hat w_x(0,t)=&q\hat{w}(0,t)-p_0 {d_m(t)} {{-p_0}\widetilde{u}(0,t), t\in \mathbb{R}_{>0}},\label{observer target system2a}\\
	{\hat w(1,t)=}&0, t\in \mathbb{R}_{>0},\label{observer target system3a}\\
	\hat w(x,0)=&\hat w_0(x),x\in(0,1), \label{observer target system4a}
\end{align}
\end{subequations}
where
\begin{align*}
{\breve \psi(x,t){:=}}&p(x){d_m(t)}-\int_{0}^{x}k(x,{z})p({z}){d_m(t)}\text{d}{z}-p_0k(x,0){d_m(t)},\\
K_p(x):=&p(x){-p_0k(x,0)}-\int_{0}^{x}k(x,{z})p({z})\text{d}{z},\\
\hat w_0(x):=&\hat u_0(x)-\int_{0}^{x}k(x,{z})\hat u_0({z})\text{d}{z}, .
\end{align*}

It is clear that system~\eqref{observer target system} has the same form as the target system~\eqref{target system} of the $u$-system. Therefore, system~\eqref{observer target system}  admits a solution $\hat w\in C^{2,1}(Q_T)\cap C\left(\overline{Q}_T\right)$  for any $T\in \mathbb{R}_{>0}$. Moreover, by virtue of Lemma~\ref{lem.1'}, the transformation \eqref{observer transfer} is invertible. Thus, it can be proved the following Lemma, whose proof is standard and hence, is omitted.
\begin{lemma}\label{observer inverse transfer}
The inverse of the transformation~\eqref{observer transfer} is given by
\begin{align}\label{observer inverse transfer conclusion}
	\hat u {(x,t)}:=\hat w {(x,t)}+\int_0^x l(x,{z}) \hat w({z},t) \mathrm{d} {z},
\end{align}
where $l$ is determined by \eqref{lem.1'}.
Moreover, the observer system~\eqref{observer system} is equivalent to its target system~\eqref{observer target system}.
\end{lemma}

In addition, for the target system~\eqref{observer target system}, we  can prove the following robust stability property, which, along with Theorem~\ref{error main result}, plays a vital role in establishing the ISS of {the  closed-loop system under output feedback control}.

\begin{proposition}\label{ISS-observer-target} The target system~\eqref{observer target system}
admits the following estimate  for all $T \in \mathbb{R}_{> 0}$:
\begin{align*}
	\left\|\hat{w}[T]\right\|_{L^{\infty}(0,1)}
	\leq&
c_4\left(e^{-{\sigma T}}\left\|\hat w_0\right\|_{L^\infty(0,1)} +e^{-\sigma T} \left\|\widetilde u_0\right\|_{L^\infty(0,1)} + \left\| f\right\|_{L^{\infty}((0, T);L^\infty(0,1))}  +\left\| d_0\right\|_{L^{\infty}(0, T)}+ \left\| d_1\right\|_{L^{\infty}(0, T)}
\right.\notag\\
&\left.+{\left\| d_m\right\|_{L^{\infty}(0, T)}} \right),
\end{align*}
where $\sigma$ is an arbitray constant satisfying $\sigma\in (0,\underline{c})$,  and $c_4$  is a positive constant depending only on $ \underline{c}$,  $\sigma$, $p_{0}$, $q$, $\max _{(x,{z})\in D}|k(x, {z})|$, $\max_{x\in[0,1]}|p(x)|$, $\max _{(x,{z})\in D}|m(x, {z})|$,  and $\max _{(x,{z})\in D}|n(x, {z})|$.
\end{proposition}

{The proof of Proposition~\ref{ISS-observer-target} is given in Section~\ref{Sec: ISS-assessment-observer-target}.}

Now, by virtue of Lemma~\ref{observer inverse transfer}, Proposition~\ref{ISS-observer-target}, and applying Theorem~\ref{error main result}, we obtain the following theorem, which {is the} third main result,  indicating the ISS of the  closed-loop system {under  output feedback control}.

\begin{theorem}\label{ISS-output-feedback-closedloop} Under the output feedback control law \eqref{observer control law},  {the   closed-loop}
{$(u,\widetilde{u})$-system}
is ISS in the  $L^\infty$-norm,
{having the following estimate  for all $T \in \mathbb{R}_{> 0}$:}
{\begin{align*}
		\left\|{u}[T]\right\|_{L^{\infty}(0,1)}+\left\|\widetilde{u}[T]\right\|_{L^{\infty}(0,1)}\leq&
		c_5\left(e^{-{\sigma  T}}\left\|u_0\right\|_{L^\infty(0,1)} +e^{-\sigma T} \left\|\widetilde u_0\right\|_{L^\infty(0,1)} + \left\|f\right\|_{L^{\infty}((0, T);L^\infty(0,1))}  +\left\| d_0\right\|_{L^{\infty}(0, T)} \right.\notag\\
		&\left.+ \left\| d_1\right\|_{L^{\infty}(0, T)}+ {\left\| d_m\right\|_{L^{\infty}(0, T)}} \right),
	\end{align*}
	where  $\sigma$ is an arbitrary constant satisfying $\sigma\in (0,\underline{c})$, and $c_5$ is a positive constant depended only on $q,  \underline{c}, \sigma, p_{0}$, $\max _{x\in [0,1]}|p(x)|$, $\max_{(x,{z})\in D}|k(x,{z})|$, $\max_{(x,{z})\in D}|l(x,{z})|$, $\max _{(x,{z})\in D}|m(x, {z})|$, and $\max_{(x,{z})\in D}|n(x,{z})|$.}
	\end{theorem}
	
	{The proof of Theorem~\ref{ISS-output-feedback-closedloop} is standard and hence, is omitted.}

\subsection{ISS Assessment of the  Closed-loop System}\label{Sec: ISS-assessment-observer-target}


\begin{pf-ISS-observer-target}  We proceed with the proof in {four} steps.

\textbf{Step 1}:  {define truncation functions.}
{Let $g$ and $G$ be defined by~\eqref{gG}.}
In addition to~\eqref{properties}, $G$ admits the following properties:
\begin{subequations}\label{gG-3}
	\begin{align}
		G(\theta ) =&\frac{1}{r+1} g(\theta ) \theta, \forall\theta  \in \mathbb{R} \label{gG3},\\
		G(\theta+\rho) \leq&2^r(G(\theta)+G(\rho)),\forall\theta,\rho\in\mathbb R.\label{G3}
	\end{align}
\end{subequations}

For any ${\sigma}\in (0,\underline {c})$ and any $T\in\mathbb{R}_{>0}$, let
\begin{align}\label{def-D}
	\mathcal{D}:=\max & \left\{\frac{e^{\sigma T}p_0}{q}\left\|{d_m}\right\|_{L^{\infty}(0, T)},\frac{e^{\sigma T}}{\underline{c}-{\sigma}}\left\|\breve{\psi}\right\|_{L^{\infty}((0,T);L^\infty(0,1))}\right\}.
\end{align}

\textbf{Step 2}: {estimate $\int_{0}^{1}G\left(e^{\sigma T}\hat w(x,T)-\mathcal{D}\right)\text{d}x$.}
Let
\begin{align*}
	\hat {v}(x, t):=e^{{\sigma} t} \hat{w}(x, t),\hat{v}_0(x):=\hat{w}_0(x),\acute{d}_0(t):= -p_0e^{{\sigma} t} {d_m(t)},
	{\acute{\psi}(x, t):= e^{{\sigma} t} \breve \psi(x, t) },
\end{align*}
{where $\hat w$ is the observer target system {(see}~\eqref{observer target system})} .
%

%
By integrating by parts and applying the definitions of $g$ and $G$, we get
\begin{align}\label{wupper}
	&\frac{\text{d}}{\text{d} t} \int_0^1 G(\hat{v}-\mathcal{D}) \text{d} x \notag\\
	=& g({\hat v}(1, t)-\mathcal{D}) {\hat v}_x(1, t)-g({\hat v}(0, t)-\mathcal{D})\left(q{\hat v}(0,t)+\acute{d} _0(t){-p_0e^{\sigma t}}\widetilde u(0,t)\right) -\int_0^1 g^{\prime}(\hat {v}-\mathcal{D}) {\hat v}_x^2 \text{d} x \notag\\
	&+\int_0^1 g({\hat v}-\mathcal{D}) (\sigma-c( t))(\hat{v}-\mathcal{D}) \text{d} x+\int_0^1 g({\hat v}-\mathcal{D}) \Big((\sigma-c( t))\mathcal{D}+\acute{\psi} \Big)+\int_0^1 g({\hat v}-\mathcal{D}) e^{\sigma t}K_p(x)\widetilde u(0,t) \text{d} x\notag\\
	=& I_1+I_2+I_3+I_4+I_5+I_{6},  \forall t\in (0,T),
\end{align}
where \begin{align*}
	{I_1}:=&g\left({\hat v}(1, t)-\mathcal{D}\right) {\hat v}_x(1, t),\\
	{I_2}:=&-g({\hat v}(0, t)-\mathcal{D})(q\hat v(0, t)+\acute{d}_0(t){-p_0e^{\sigma t}}\widetilde u(0,t)),\\
	{I_3}:=&-\int_0^1 g^{\prime}\left({\hat v}(x,t)-\mathcal{D}\right) {\hat v}_x^2(x,t)  \text{d} x,\\
	{I_4}:=&\int_0^1 g({\hat v}(x,t)-\mathcal{D}) (\sigma-c( t))(\hat{v}(x,t)-\mathcal{D}) \text{d} x,\\
	{I_5}:=&\int_0^1 g({\hat v}(x,t)-\mathcal{D}) \Big((\sigma-c( t))\mathcal{D}+\acute{\psi}(x,t) \Big) \text{d} x,\\
	{I_{6}}:=&\int_0^1 g({\hat v}(x,t)-\mathcal{D})K_p(x)e^{\sigma t}\widetilde u(0,t) \text{d} x.
\end{align*}
It is easy to see that {$I_1=0$}, {$I_3\leq 0$}, and {$I_5\leq 0$} for all $t\in(0,T)$. In addition, by the definitions of $\mathcal{D}$ and $ \acute{d}_0$, the Young's inequality with $\varepsilon_1\in(0,q)$, and \eqref{gG3}, we deduce  that
\begin{align}\label{i6}
	{I_2}
	\leq&-q(r+1) G(\hat{v}(0, t)-\mathcal{D})+r\varepsilon_1 G(\hat{v}(0, t)-\mathcal{D})+\varepsilon_1^{-r} G\left({|p_0|}e^{\sigma t} |\widetilde{u}(0, t)|\right) \notag\\
	\leq& \frac{\varepsilon_1^{-r}}{r+1} {|p_0|}^{r+1} e^{(r+1) \sigma t}|\widetilde{u}(0, t)|^{r+1},\forall t\in (0,T).
\end{align}

{For $I_4$}, by $\underline{c}-{\sigma}>0$ and the properties of $g$ and $ G$ (see \eqref{gG3}), we obtain
\begin{align}\label{i8}
	{I_4}
	\leq-\left(\underline{c}-{\sigma}\right)(r+1) \int_0^1 G(\hat{v}-\mathcal{D}) \text{d} x,\forall t\in (0,T).
\end{align}

For {$I_{6}$}, applying the Young's inequality, \eqref{gG3}, and the definition of $g$ (see \eqref{gG}), we get for all $t\in (0,T)$:
\begin{align}\label{i10}
	{I_{6}}
	\leq & r\varepsilon_2\int_0^1 G(\hat{v}-\mathcal{D}) \text{d} x
	+\frac{\varepsilon_2^{-r} }{r+1}\overline{K}_p^{r+1}e^{ (r+1)\sigma t}|\widetilde u(0,t)|^{r+1},
\end{align}
where  $\varepsilon_2$ is a positive constant that will  be determined later and $\overline{K}_p:= \max_{x\in[0,1]}|K_p(x)| $.

Let $\varepsilon_2\in (0, \underline{c}-{\sigma})$ and $\varepsilon_3:=\min\{\varepsilon_1,\varepsilon_2\}$.
Since {$I_i\leq0~(i=1,3,5)$}, by \eqref{wupper}, \eqref{i6}, \eqref{i8}, and~\eqref{i10}, we get for all $t\in (0,T)$:
\begin{align*}
	 \frac{\text{d}}{\text{d} t} \int_0^1 G(e^{\sigma t}\hat{w}-\mathcal{D}) \text{d} x=& \frac{\text{d}}{\text{d} t} \int_0^1 G(\hat{v}-\mathcal{D}) \text{d} x \\
	\leq &-r_0 \int_0^1 G\left(e^{\sigma t}\hat{w}-\mathcal{D}\right) \text{d} x +\frac{ \varepsilon_3^{-r} K_0^{r+1} }{r+1}       e^{(r+1) \sigma t}
	|\widetilde{u}(0, t)|^{r+1},
\end{align*}
where $r_0:=r_0(\varepsilon_2,r):=\left(\underline{c}-{\sigma}\right)(r+1)-r\varepsilon_2>0$ and $K_0:= {|p_0|}+ \overline{K}_p  $.

By using the Gronwall's inequality,
we deduce that
\begin{align*}
	 \int_0^1 G(e^{\sigma T}\hat{w}(x,T)-\mathcal{D}) \text{d} x
	\leq e^{-r_0 T} \int_0^1 G\left(\hat{w}_0(x)-\mathcal{D}\right) \text{d} x+\frac{ \varepsilon_3^{-r} K_0^{r+1} }{r+1} \int_0^Te^{-r_0(T-s)+(r+1)\sigma s}  |\widetilde{u}(0,s)|^{r+1} \text{d}s.
\end{align*}
Note that  the spatial $L^{\infty}$-estimate of $\widetilde{u} $ (see {Theorem~\ref{error main result}}) ensures that
\begin{align*}
	|\widetilde{u}(0, T)|^{r+1}
	& \leq 2^r{C_1}^{r+1} \left(e^{-\sigma(r+1) T}\left\|\widetilde{u}_0\right\|_{L^\infty(0,1)}^{r+1}+B_1^{r+1}\right) ,
\end{align*}
where
\begin{align}\label{def-b1}
	{B_1:=}&{\left\|{d_m}\right\|_{L^{\infty}(0, T)}+}\left\|{d}_0\right\|_{L^{\infty}(0, T)}+\left\|{d}_1\right\|_{L^{\infty}(0, T)}+\left\|{f}\right\|_{L^{\infty}((0, T) ; L^{\infty}(0,1))},
\end{align}
and {$C_1$}  {is a positive constant depending only on $ \underline{c},  \sigma, p_{0}, {q}$,   {$\max_{x\in[0,1]}|p(x)|$,}  $\max _{(x,{z})\in D}|m(x, {z})|$, and {$\max _{(x,{z})\in D}|n(x, {z})|$}.} Then, we have
\begin{align}\label{hatw-D}
	 \int_0^1 G\left(e^{\sigma T}\hat{w}(x,T)-\mathcal{D}\right) \text{d} x
	\leq& e^{-r_0 T} \int_0^1 G\left(|\hat{w}_0(x)|-\mathcal{D}\right) \text{d} x+\frac{{C_2}}{r_0(r+1)}\left\|\widetilde{u}_0\right\|_{L^{\infty}(0,1)}^{r+1} \notag\\
	&+\frac{{C_2}e^{(r+1) \sigma T}}{(r_0+(r+1) \sigma)(r+1)} B_1^{r+1},
\end{align}
where ${C_2:=C_2 (\varepsilon_3,r):=2^rC_1^{r+1}}\varepsilon_3^{-r}K_0^{r+1}$. 

\textbf{Step 3:}  estimate $\int_{0}^{1}G(-e^{\sigma T}\hat w(x,T)-\mathcal{D})\text{d}x$.
Analogous to the derivations in Step~2, by calculating $\frac{\text{d}}{\text{d}t}\int_0^1G\left(-e^{\sigma T}\hat{w}(x, T)-{{\mathcal{D}}}\right) \text{d} x$ and using the Gronwall's inequality, we can obtain
\begin{align}\label{-hatw-D}
	\int_0^1G\left(-e^{\sigma T}\hat{w}(x, T)-{{\mathcal{D}}}\right) \text{d} x\leq& e^{-r_0 T} \int_0^1 G\left(|\hat{w}_0(x)|-\mathcal{D}\right) \text{d} x+{\frac{C_2}{r_0(r+1)}}\left\|\widetilde{u}_0\right\|_{L^{\infty}(0,1)}^{r+1}  \notag\\
	&+\frac{{C_2}e^{(r+1) \sigma T}}{(r_0+(r+1) \sigma)(r+1)} B_1^{r+1}.
\end{align}

\textbf{Step 4:} establish the estimate of $\|\hat w[T]\|_{L^\infty(0,1)}$.  By the   properties  of $g$ and $G$ (see \eqref{gG-3}),  \eqref{hatw-D}, and~\eqref{-hatw-D}, we get
\begin{align*}
	\frac{e^{(r+1) \sigma T}}{r+1} \int_0^1|\hat{w}(x, T)|^{r+1} \text{d} x
	\leq& 2^r\left(\int_0^1 G\left(e^{\sigma T}|\hat{w}(x, T)|-{\mathcal{D}}\right)\text{d} x+\int_0^1 G({\mathcal{D}}) \text{d} x\right) \\
	\leq& \frac{2^r}{r+1} e^{-r_0 T}\int_0^1 |\hat{w}_0(x)|^{r+1} \text{d} x  +\frac{2^r{C_2}}{r_0(r+1)}\left\|\widetilde{u}_0\right\|_{L^{\infty}(0,1)}^{r+1} \notag\\
	& +\frac{2^re^{(r+1) \sigma T}{C_2}}{(r_0+(r+1) \sigma)(r+1)} B_1^{r+1} +\frac{2^r}{r+1}  {\mathcal{D}}^{r+1} ,
\end{align*}
which 
 implies that
\begin{align}
	\|\hat{w}[T]\|_{L^{r+1}(0,1)}
	\leq& 2^{\frac{r}{r+1}} e^{-\left(\frac{r_0}{r+1}+\sigma\right) T}\left\|\hat{w}_0\right\|_{L^{r+1}(0,1)}+\left(\frac{2^r{C_2}}{r_0}\right)^\frac{1}{r+1}e^{-\sigma T}\left\|\widetilde{u}_0\right\|_{L^{\infty}(0,1)}\notag\\ &+\left(\frac{2^r{C_2}}{r_0+(r+1) \sigma}\right)^{\frac{1}{r+1}}B_1 +2^{\frac{r}{r+1}}e^{-\sigma T}{\mathcal{D}}.\label{equ.45}
\end{align}
	Letting $r\to +\infty$ in \eqref{equ.45}, we obtain
	\begin{align}\label{uuuu}
		 \|\hat{w}[T]\|_{L^{\infty}(0,1)}
		\leq 2 e^{-\left(\underline{c}-\varepsilon_2\right){T}}\left\|\hat{w}_0\right\|_{L^{\infty}(0,1)}+{C_3}e^{-\sigma T}\left\|\widetilde{u}_0\right\|_{L^{\infty}(0,1)}
		+{C_3} B_1+2 e^{-\sigma T} {\mathcal{D}},
	\end{align}
	where ${C_3}:=4{C_1}K_0\varepsilon_3^{-1}$ and {$\underline{c}-\varepsilon_2>\sigma>0$}. Then, in view of the definitions of $B_1$ (see~\eqref{def-b1}) and $\mathcal{D}$ (see~\eqref{def-D}), we conclude that Proposition~\ref{ISS-observer-target} holds true.
\end{pf-ISS-observer-target}

\section{Numerical Results}\label{sec:numerical results}
In this section, we present some numerical results on the ISS in the $L^\infty$-norm for the considered systems  with different disturbances under the state feedback control law~\eqref{control law}, and the observer-based output feedback control law~\eqref{observer control law}, respectively.

In simulations, we always set
\begin{align*}
	q =&1,\\
	\lambda(t) = &\left\{
	\begin{aligned}
		&1.2\pi^2\left( \sin^2 5t +1\right), 0\le t\le1,\\
		&1.2\pi^2\left( \sin^2 5t +e^{-t}-e^{-1}+1\right), t > 1.
	\end{aligned}
	\right.
\end{align*}
In \eqref{c0} and \eqref{specific-k}, we set $c_0=\frac{13}{5}\pi^2$.  Let {$p_0=\frac{6}{5}\pi^2$}. {The gain $p(x)$   is obtained by using the successive approximations.} 

Moreover, we set
	\begin{align*}
		{u_0}(x)=&-\frac{1}{6}\left(5x-\frac{1}{4}\right)(2-x)\left(3x^2-1\right),\\
		{\hat {u}_0}(x)=& 0,\\
		{\widetilde{u}_0}(x)=& u_0(x)-\hat u_0(x), \\
		f(x,t)=& \frac{A}{6} \sin\left(60t+x\right),\\
		d_0(t) =&  \frac{A_0}{5} \left( \sqrt{t}  e^{-t} +2\sin(25t)\right) ,\\ 
		d_1(t)=&   \frac{2}{5}A_1 \sin(25t),\\
  {d_m(t)}= &\frac{A_2}{44}\sin(40t),
	\end{align*}
	where $A\in\{0,2\}$ and $A_0, A_1, A_2\in \{0,1,3\}$ are used to describe
	the amplitude of disturbances.

\begin{figure}[htbp]
	\centering
	\begin{minipage}{1\linewidth}	
		\subfigure[Evolution of  $u$]{
			\includegraphics[scale=0.38]{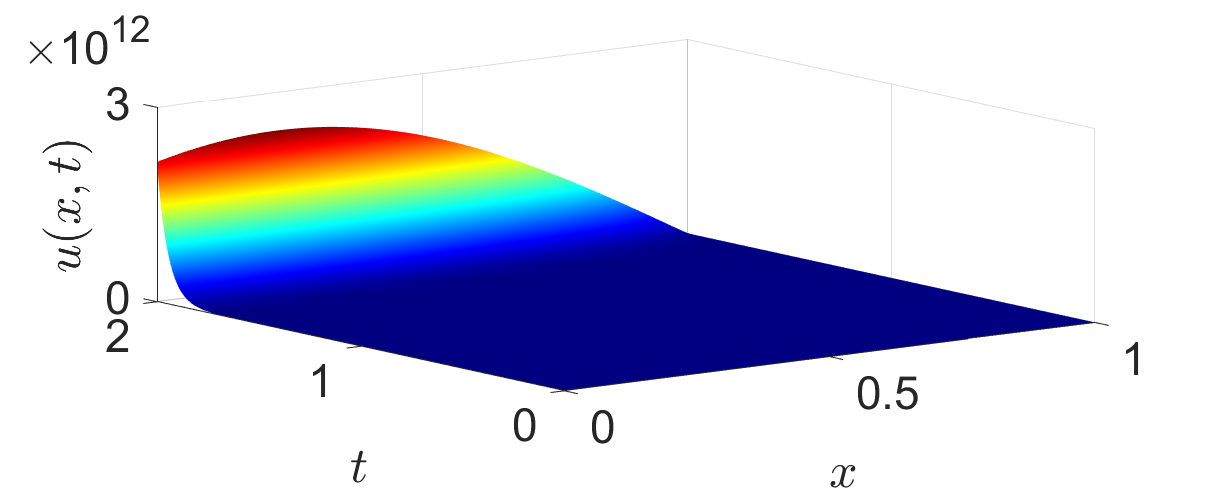}	
		}\noindent
		\subfigure[ {Evolution of  $\|u[t]\|_{L^{\infty}(0,1)}$}]{
			\includegraphics[scale=0.38]{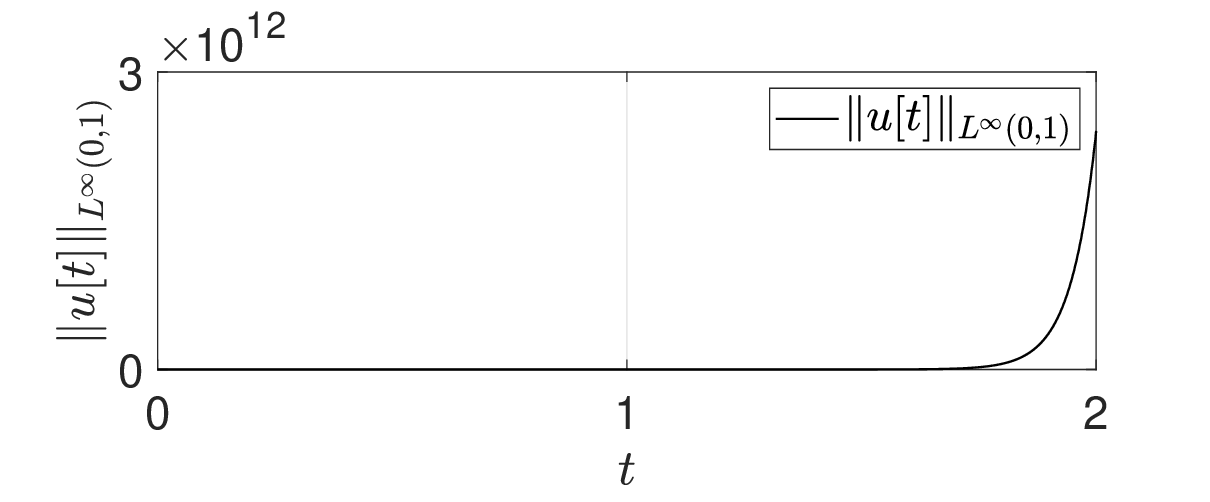}
		}
	\end{minipage}
	\centering
	\caption{Open-loop response of   the $u$-system~\eqref{original system} without disturbances}
	\label{fig1}
\end{figure}

\begin{figure}[htbp]
	\begin{center}
		\begin{minipage}{1\linewidth}	
			\subfigure[Evolution of   $u$ when   $A=A_0=A_1=0$]{
				\includegraphics[scale=0.38]{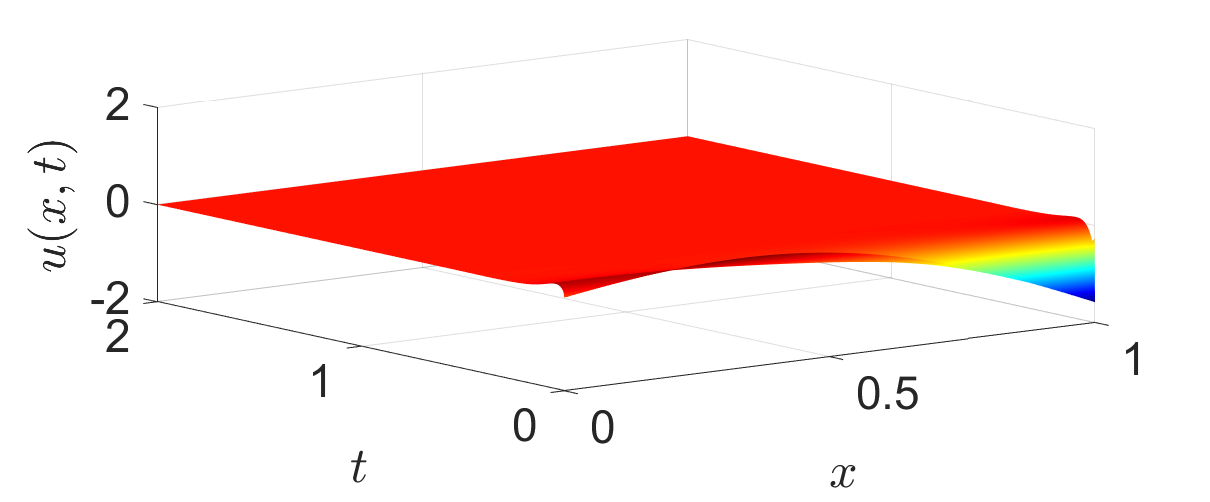}	
			}\noindent
			\subfigure[ Evolution of   $u$ when   $A=2,A_0=A_1=1$]{
				\includegraphics[scale=0.38]{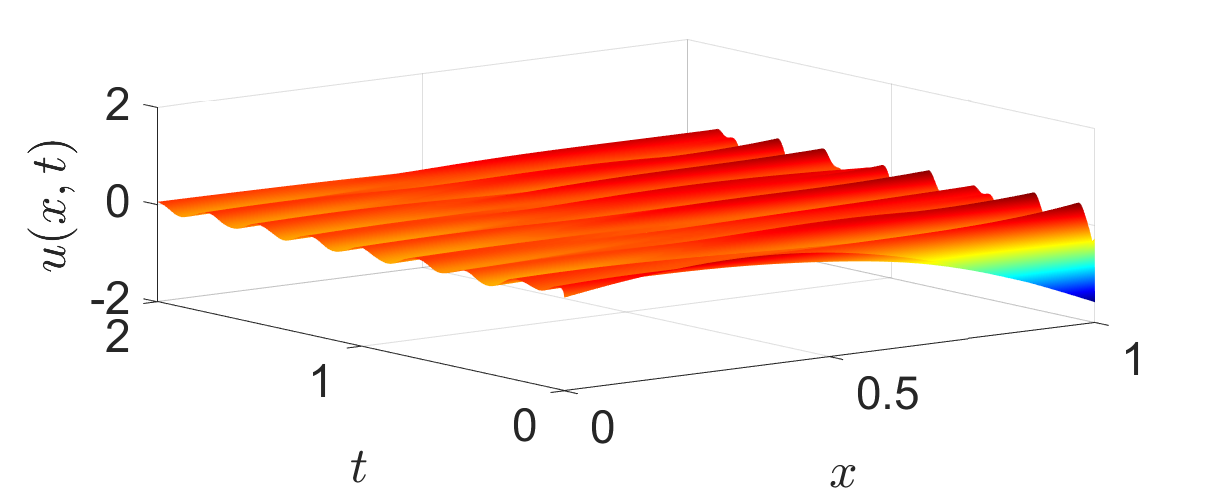}
			}
		\end{minipage}
		\begin{minipage}{1\linewidth }
			\subfigure[ Evolution of   $u$ when   $A=2,A_0=A_1=3$]{
				\includegraphics[scale=0.38]{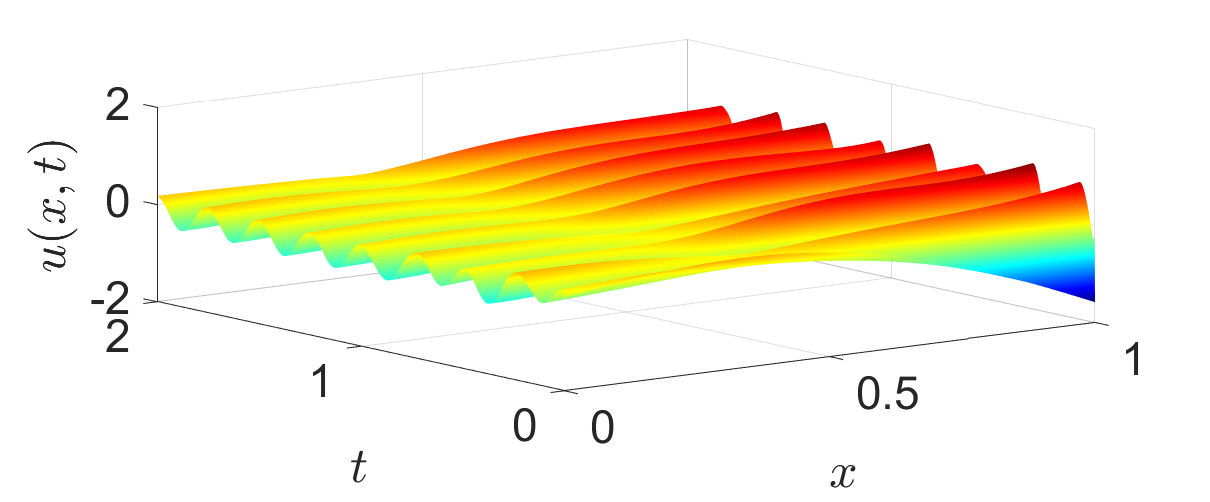}
			}\noindent
			\subfigure[{Evolution of  $\|u[t]\|_{L^\infty(0,1)}$ for different disturbances}]{
				\includegraphics[scale=0.38]{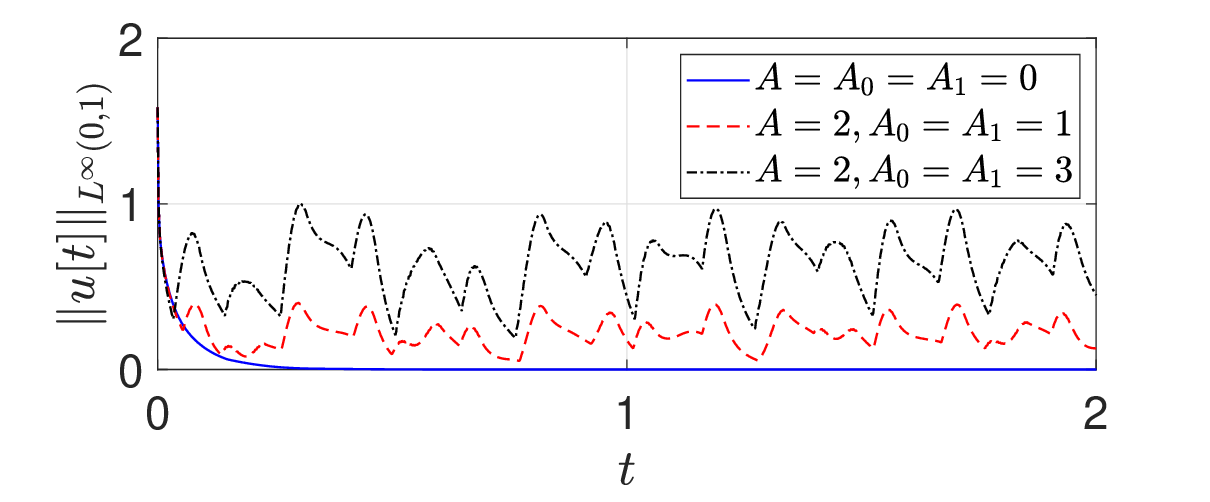}
			}
		\end{minipage}
		\caption{Closed-loop responses of the $u$-system~\eqref{original system} with different disturbances under the state feedback control law~\eqref{control law} }
		\label{fig2}
	\end{center}
\end{figure}

\begin{figure}[htbp]
	\centering
	\begin{minipage}{1\linewidth}	
		\subfigure[Evolution of   $\widetilde u$ when   $A=A_0=A_1=A_2=0$]{
			\includegraphics[scale=0.38]{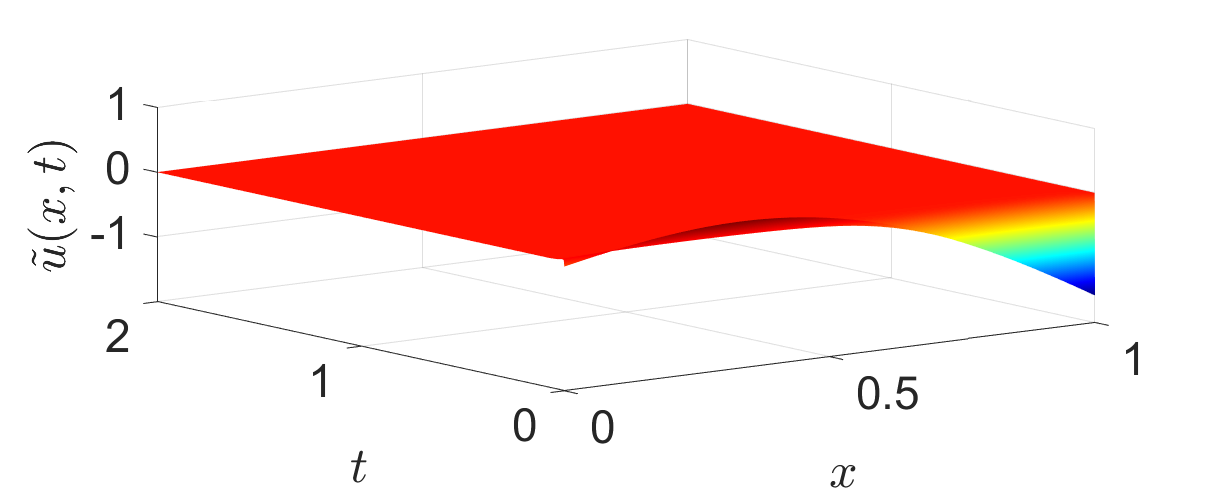}	
		}\noindent
		\subfigure[ Evolution of   $\widetilde u$ when   $A=2,A_0=A_1=A_2=1$]{
			\includegraphics[scale=0.38]{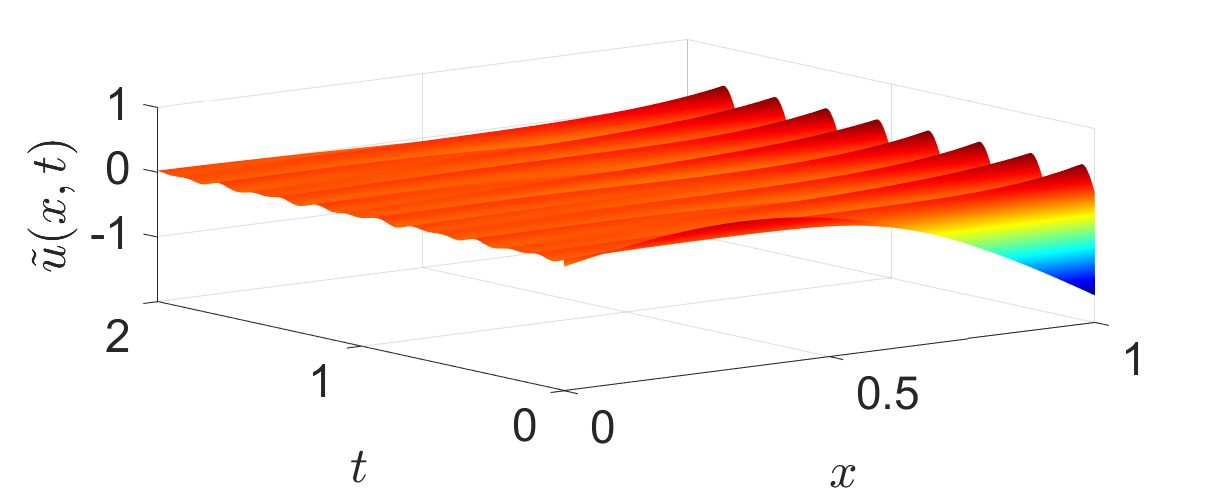}
		}
	\end{minipage}
	\begin{minipage}{1\linewidth }
		\subfigure[ Evolution of   $\widetilde u$ when   $A=2,A_0=A_1=A_2=3$]{
			\includegraphics[scale=0.38]{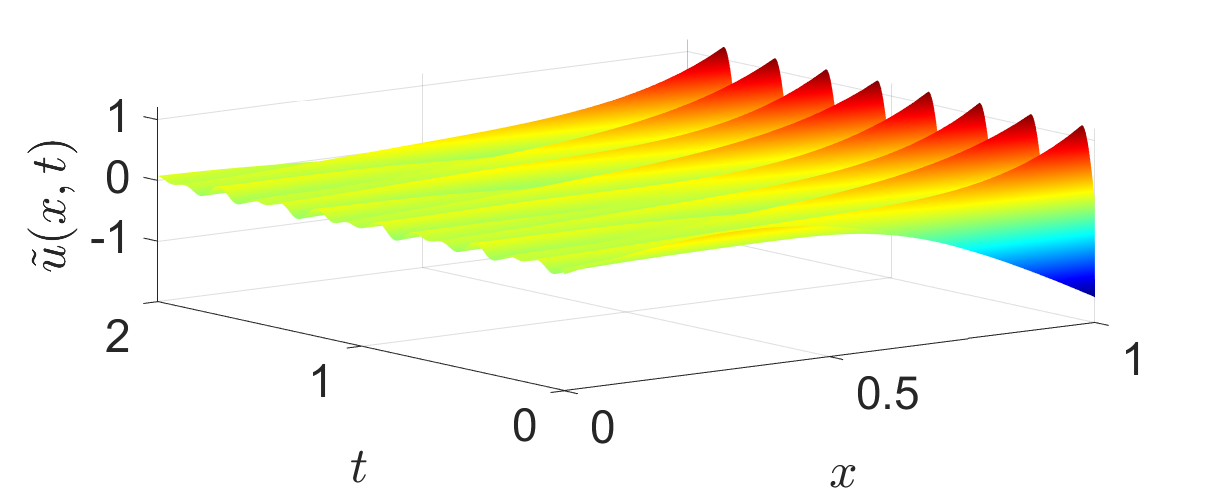}
		}\noindent
		\subfigure[{Evolution of  $\|\widetilde u[t]\|_{L^\infty(0,1)}$ for different disturbances}]{
			\includegraphics[scale=0.38]{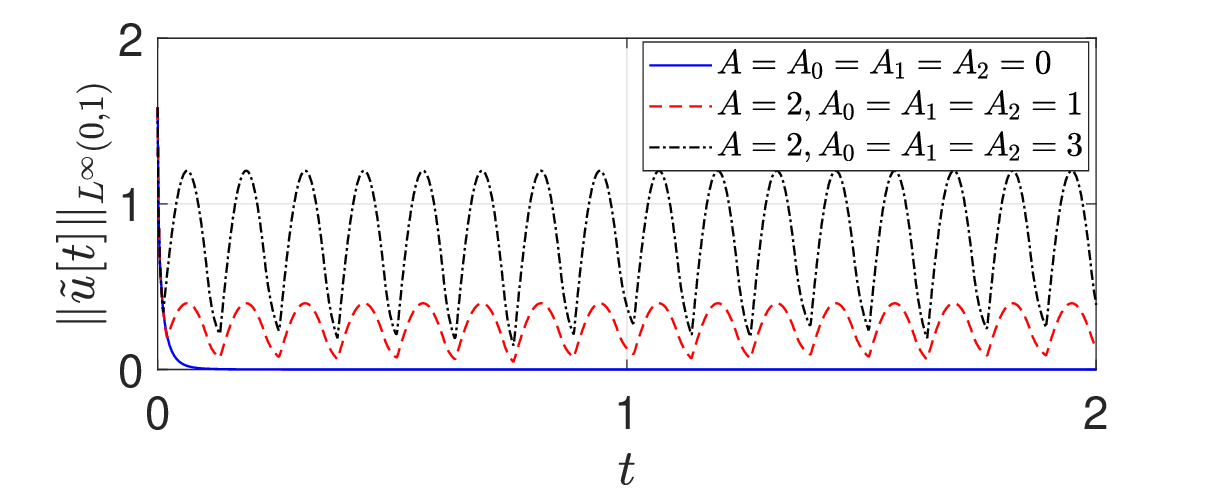}
		}
	\end{minipage}
	\caption{Evolution of   $\widetilde u$ {and $\|\widetilde u[t]\|_{L^\infty(0,1)}$ }for the error $\widetilde u$-system~\eqref{error system} with  different disturbances }
	\label{fig3}
\end{figure}

Figure~\ref{fig1} shows that the disturbance-free $u$-system~\eqref{original system} is unstable {in open loop}, while Fig.~\ref{fig2}(a), (b), and (c) show that, under the state feedback control law~\eqref{control law},  the  state of the disturbance-free system~\eqref{original system} tends to the origin and    the    states of system~\eqref{original system} with different  disturbances remain  bounded. Moreover, when  amplitudes of disturbances  decrease, the  amplitudes of states  decrease. This reflects well the ISS property of the $u$-system under the state feedback control law~\eqref{control law} as shown in  Fig.~\ref{fig2}(d), where the evolution of $\|u[t]\|_{L^\infty(0,1)}$ {w.r.t.}  different {disturbances} are presented.
{Similarly}, {Fig.}~\ref{fig3} reflects   the ISS property of {the $\widetilde{u}$}-system~\eqref{error system}.

\begin{figure}[htbp]
	\centering
	\begin{minipage}{1\linewidth}
		\subfigure[    Evolution of   $u$ when   $A=A_0=A_1=A_2=0$]{
			\includegraphics[scale=0.38]{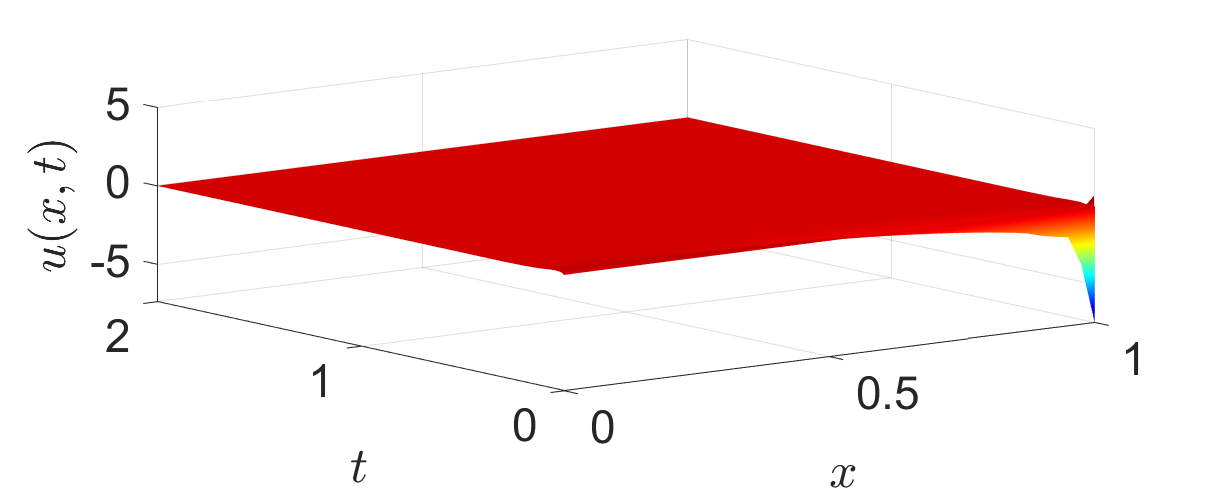}	
		}\noindent
		\subfigure[ Evolution of   $u$ when   $A=2,A_0=A_1=A_2=1$]{
			\includegraphics[scale=0.38]{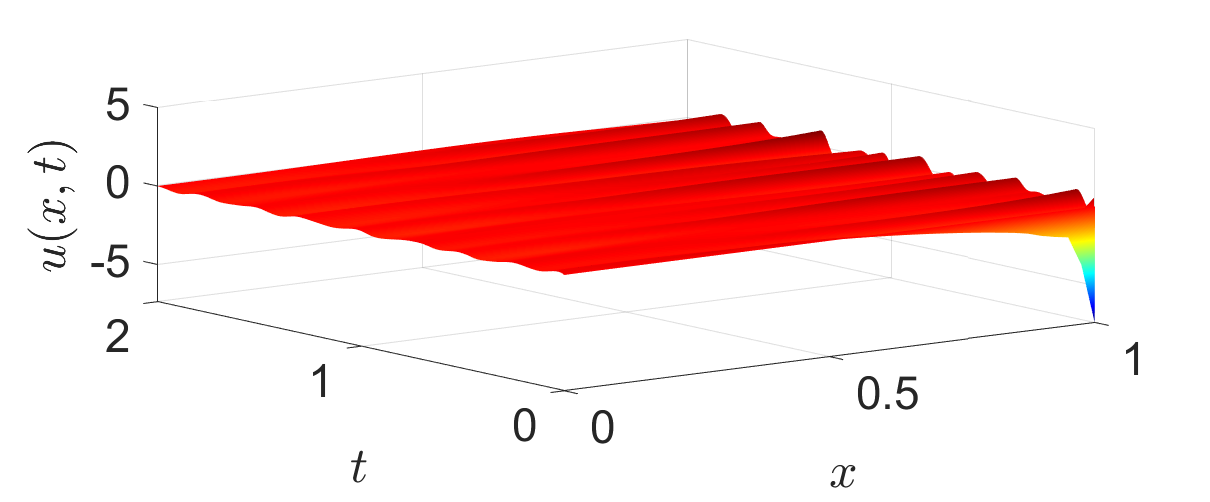}
		}
	\end{minipage}
		\centering
		\includegraphics[scale=0.38]{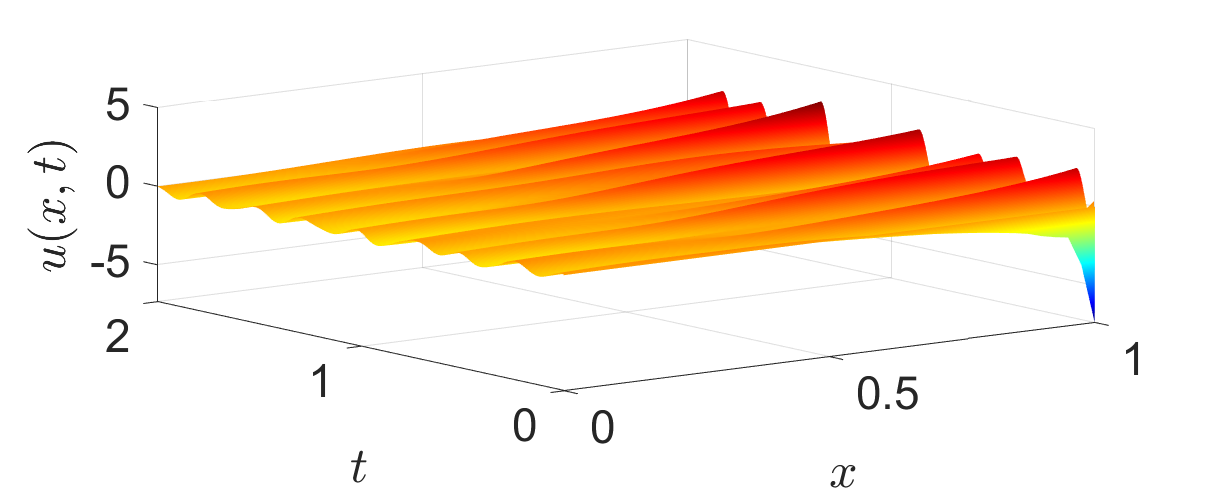}\\
	\small{(c)  Evolution of   $u$ when   $A=2, A_0=A_1=A_2=3$}
	\caption{Closed-loop responses of the $u$-system~\eqref{original system} with different disturbances under the output feedback control law~\eqref{observer control law}}
	\label{fig4}
\end{figure}

\begin{figure}[htbp]
	 \centering
	\begin{minipage}{1\linewidth}
		\subfigure[Evolution of   $\hat u$ when   $A=A_0=A_1=A_2=0$]{
			\includegraphics[scale=0.38]{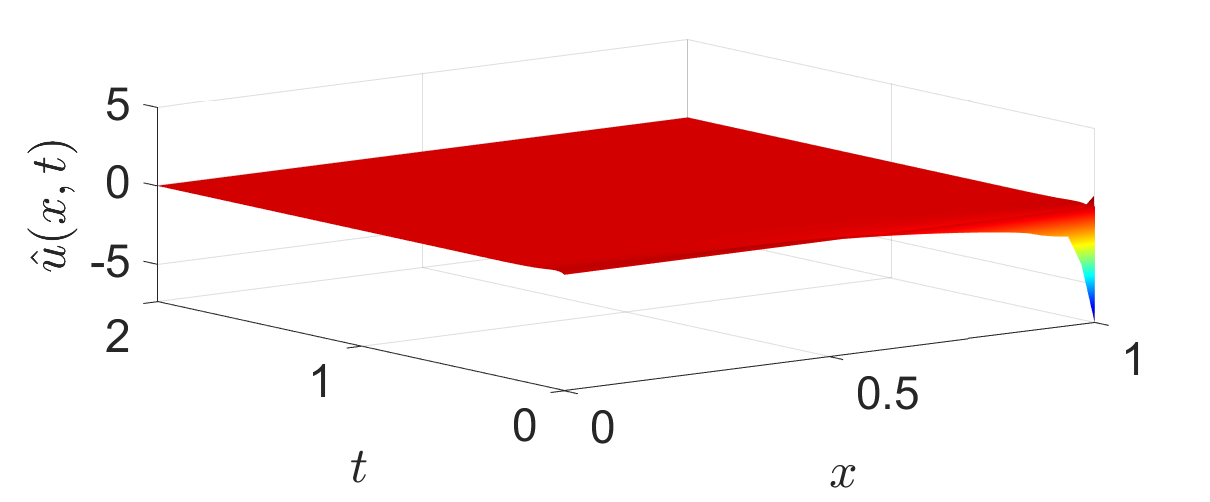}	
		}\noindent
		\subfigure[ Evolution of   $\hat u$ when   $A=2,A_0=A_1=A_2=1$]{
			\includegraphics[scale=0.38]{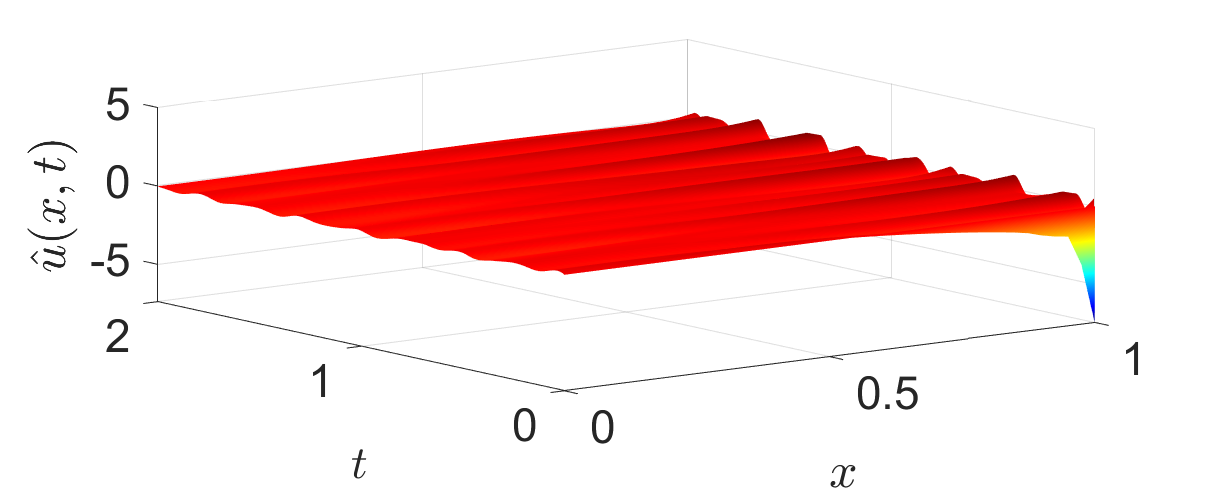}
		}
	\end{minipage}
	\begin{minipage}{1\linewidth }
		\subfigure[ Evolution of   $\hat u$ when   $A=2,A_0=A_1=A_2=3$]{
			\includegraphics[scale=0.38]{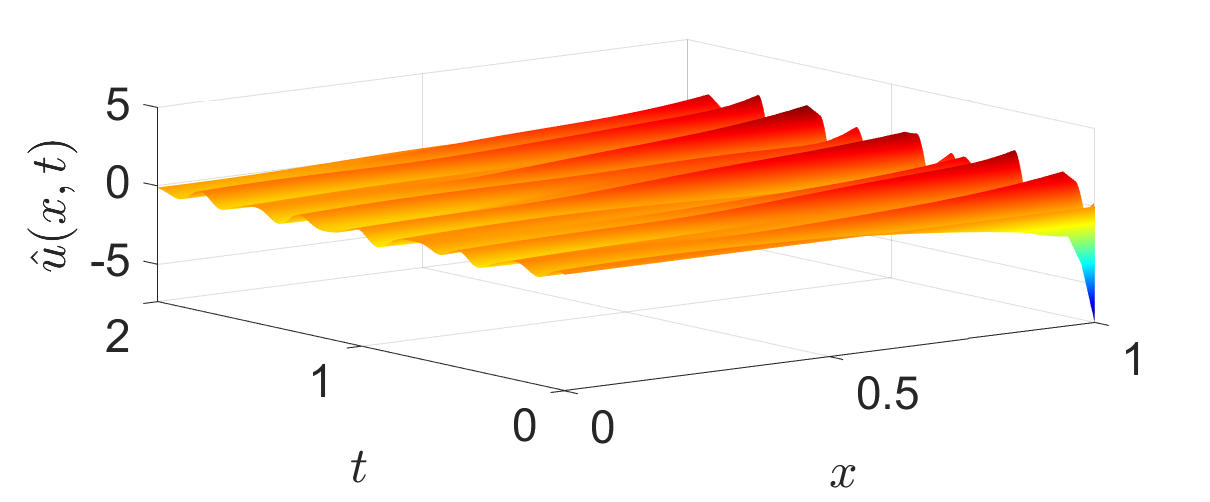}
		}\noindent
		\subfigure[{Evolution of $\|  u[t]\|_{L^\infty(0,1)}+\|\widetilde u[t]\|_{L^\infty(0,1)}$ for different disturbances}]{
			\includegraphics[scale=0.38]{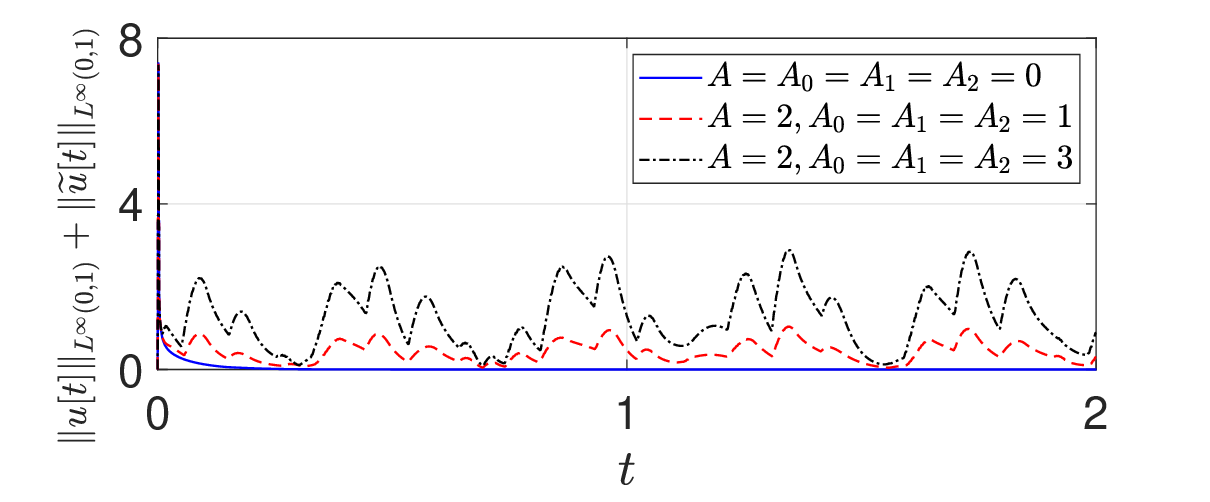}
		}
	\end{minipage}
	\caption{{Evolution of   $\hat u$ and $\|  u[t]\|_{L^\infty(0,1)}+\|\widetilde u[t]\|_{L^\infty(0,1)}$} for the observer  $\hat u$-system~\eqref{observer system} and the  $(u,\widetilde u)$-system coupled via the output feedback control law~\eqref{observer control law}}
	\label{fig5}
\end{figure}

Under the output feedback control law~\eqref{observer control law}, Fig.~\ref{fig4} indicates   the robust {stability} of the $u$-system~\eqref{original system}, and
Fig.~\ref{fig5}(a), (b), and (c) {show}   the robust {convergence} of the   observer $\hat u$-system~\eqref{observer system}, respectively. In particular, {the results shown in }Fig.~\ref{fig5}(d) {reflect} well the ISS property of {the} closed-loop $(u, \widetilde{u})$-system. It can be seen that {the} numerical results {are}  consistent with the theoretical  {analysis}, {demonstrating well} the effectiveness of the proposed {control scheme}.

\section{Conclusion}\label{conclusion}
In this paper, we designed {an} observer-based {output feedback stabilizing controller} for a class of time-varying parabolic PDEs with in-domain and mixed boundary disturbances under the framework of ISS. {Time-invariant} kernel functions    {were} used in both controller and observer design and implementation, which required much less computation efforts compared to those required in the standard backstepping control of PDEs.  {The} generalized Lyapunov method, which {allows dealing} with Dirichlet boundary disturbances, was applied to establish the ISS in the $L^{\infty}$-norm of {the considered system} {in closed loop}. To validate the effectiveness of the proposed {control scheme},  numerical simulations were conducted for systems with {different} disturbances.
It is worth mentioning that establishing the ISS in the $L^p$-norm with $p\in[2,+\infty)$ for the   considered {system in  closed loop under   output feedback control is much challenging due to the appearance of {the boundary term $\widetilde{u}(0,t)$}, which    is difficult to estimate by using the $L^p$-norm of $\hat{w}$ or $\widetilde{u} $, in the target system~\eqref{observer target system} of the observer.} This will be considered in our future work.

\end{document}